\documentclass{amsart}
\usepackage{amsmath,amssymb,amsfonts,mathtools,tikz,subcaption,hyperref,
  mathrsfs,algorithm,algpseudocode,graphicx,etoolbox,soul,csquotes,ifthen,
  subcaption,pgfplots,array,diagbox,multirow}
\usepackage[scr=esstix]{mathalpha}
\usepackage{xcolor}
\pgfplotsset{compat=1.18}
\usepgfplotslibrary{external}
\usetikzlibrary{external,arrows.meta, positioning}
\usepackage[foot]{amsaddr}
\usepackage[utf8]{inputenc} 
\usepackage[T1]{fontenc}
\usepackage[toc,page]{appendix} 
\usepackage[space,noadjust,nocompress]{cite}

\newtheorem{theorem}{Theorem}[section]
\newtheorem{defi}[theorem]{Definition}
\newtheorem{prop}[theorem]{Proposition}
\newtheorem{lem}[theorem]{Lemma}
\newtheorem{cor}[theorem]{Corollary}
\newtheorem{rem}[theorem]{Remark}

\newcommand\scalemath[2]{\scalebox{#1}{\mbox{\ensuremath{\displaystyle#2}}}}

\begin{document} 

\title[On an efficient line smoother for the p-multigrid $\gamma$-cycle]
{On an efficient line smoother for the p-multigrid $\gamma$-cycle}

\author[J. P. Lucero Lorca]{José Pablo Lucero Lorca$^{1,2}$}
\email{josepablo.lucerolorca@colostate.edu}
\address{$^1$
  Cooperative Institute for Research in the Atmosphere (CIRA),
  Colorado State University, Fort Collins, CO 80309. \\
  \url{https://pablo.world/mathematics}}
\address{$^2$
  Affiliate working on a cooperative agreement/grant stationed at NOAA/Global \\
  Systems Laboratory, Boulder, CO.}
\author[D. Rosenberg]{Duane Rosenberg$^{1,2}$}
\email{Duane.Rosenberg@colostate.edu}
\author[I. Jankov]{Isidora Jankov$^3$}
\address{$^3$
  National Oceanic and Atmospheric Administration, Global Systems
  Laboratory R/GSL 325 Broadway, Boulder, CO 80305, USA. \\
  \url{https://gsl.noaa.gov/profiles/isidora.jankov}}
\email{isidora.jankov@noaa.gov}
\author[C. McCoid]{Conor McCoid$^4$}
\address{$^4$
  McMaster University}
\email{mccoidc@mcmaster.ca}
\author[M. Gander]{Martin Gander$^5$}
\address{$^5$
  Université de Genève
}
\email{martin.gander@unige.ch}

\maketitle

\begin{abstract}
  As part of the development of a Poisson solver for the spectral
  element discretization used in the GeoFluid Object Workbench
  (GeoFLOW) code, we propose a solver for the linear system arising
  from a Gauss-Legendre-Lobatto global spectral method. We
  precondition using a $p$-multigrid $\gamma$-cycle with
  highly-vectorizable smoothers, that we refer to as \emph{line
  smoothers}. Our smoothers are restrictions of spectral and finite
  element discretizations to low-order one-dimensional problems along
  lines, that are solved by a reformulation of cyclic reduction as a
  direct multigrid method. We illustrate our method with numerical
  experiments showing the apparent boundedness of the iteration count
  for a fixed residual reduction over a range of moderately deformed
  domains, right hand sides and Dirichlet boundary conditions.
\end{abstract}

\section{Introduction}

\subsection{Application background and motivation}

Elliptic solutions can be challenging for discretizations that utilize
local basis functions (e.g., finite volume or finite difference
methods), and yet it is precisely these schemes that are usually best
suited to discretizing terrain when examining the effects of the
atmospheric boundary layer on processes aloft.  Atmospheric flows are
typically very slow compared to the speed of sound, and are thus well
approximated by the incompressible Navier Stokes equations (e.g.,
\cite{Pope2000}). In order to maintain the incompressibility
constraint an elliptic problem needs to be solved at each time step of
the numerical solution (e.g., \cite{DevilleFischerMund2002}).
However, such flows are not strictly incompressible, which justifies
the use of the compressible Navier Stokes equations in some contexts.
One such instance is the numerical solution of atmospheric flows at
low Mach number, in order to examine multi-scale properties. This
choice has the advantage of obviating the need to solve an elliptic
problem at {\it each} timestep, but instead, enables us to do 
so sporadically.

Under atmospheric conditions, the conservation laws admit solutions
that often degenerate into turbulence by way of nonlinear interactions
that transfer energy, moisture, and aerosols irreversibly from large
to small scales.  The large number of degrees of freedom required to
capture accurately this transfer and the associated mixing mechanisms
lead to stochastic behavior in atmospheric flow simulations that
develops organically from the conservation laws alone. The GeoFluid
Object Workbench framework, called GeoFLOW \cite{Rosenberg2023}, was
developed principally to focus on effects of numerical
characteristics, such as spatial and temporal truncation and
dissipation processes, on these dynamical statistical properties.

GeoFLOW partitions the domain into a union of finite elements of
arbitrary order, in which functions are represented as expansions in
terms of a tensor product of one dimensional Lagrange interpolating
polynomials that serve as the (element-) local basis functions. The
code enables specification of terrain profiles, the mesh quality is
improved by setting a discretized Poisson equation problem on it (see
\cite{Rosenberg2023} for details) and the corresponding linear system
is solved using an iterative Krylov method. The influence of terrain
on the smooth deformation of the interior grid is dictated by the
fine-scale character of the terrain profile and it can be quite expensive
to compute without preconditioning. Nevertheless, this terrain
computation is done only at initialization, and thus, it needs only be
efficient {\it enough} to be folded into the overall start-up
cost. While terrain undeniably serves to motivate this work, it does
not govern completely the need for a highly efficient preconditioner.

For our current GeoFLOW applications the primary factor motivating
the need for an efficient preconditioner is the computation of 
specific diagnostics of the solution, as opposed to the solution 
itself. As mentioned, the atmosphere is slightly compressible, and we 
must provide the ability to distinguish between both compressible 
{\it and} incompressible modes in any attempt to unravel their effects. 
For this, we will use a Helmholtz decomposition (e.g., \cite{Passot1987}), 
which also requires a solution to a Poisson problem discretized potentially 
on a grid that is now deformed by terrain. While these diagnostics are 
computed typically at cadences of $\mathcal{O}(100-1000)$ timesteps, this 
cadence is user--specifiable, and, depending on the phenomenology under 
consideration, {\it may} be required at still higher values to capture 
high frequency signatures of these modes during a run. It is for these 
high-cadence decompositions that we mainly seek an efficient 
Poisson solver.

\subsection{Computational needs and approach}
GeoFLOW uses a spectral element method  (SEM)
consisting of a spectral Gauss-Legendre-Lobatto (GLL) spatial
discretization in each element, and an inter-element Dirichlet
coupling. As noted in the previous section, our applications demand 
that an efficient Poisson solver be employed.

Achieving efficiency implies a fast time-to-solution, but for codes
that are designed for a large degree-of-freedom count, it also implies
\emph{scalability} in the sense that when adding degrees of freedom,
the computational complexity of the solver as a whole will not {\bf
  grow} faster than the simple application of the discretized
Laplacian operator. Our choice to achieve such scalability for
elliptic problems is a multigrid--preconditioned (see
\cite{Hackbusch1985}) GMRES solver. A multigrid preconditioner
consists of a set of \emph{smoothers}, \emph{coarse spaces},
\emph{restriction} and \emph{prolongation} operators acting on
residuals, tailored to the discretized operator being
preconditioned.

The ultimate objective is obtaining a preconditioned solver that, by
using modern matrix-free, vectorized techniques, is able to leverage
parallelization to accelerate the solution, while at the same time
keeping memory footprint controlled. It is also desirable that the
method be flexible, so memory footprint and
parallelization can be tuned to different architectures.

\subsection{Context and focus of this manuscript}
This manuscript describes the first step in the development of the
preconditioned Poisson solver, the capacity to solve a spectral
problem {\bf on a single 2D element} for relatively high polynomial
degrees.

Using an efficient tridiagonal solver that we describe in detail, we
characterize two different smoothers based on the incomplete
factorization along lines in the $x$ and $y$ directions. One smoother
factorizes the original system matrix arising from the spectral
method, and the other factorizes a bilinear finite element matrix on
the GLL mesh.

Our objective is to leverage the tridiagonal solver, exploring the
computational complexity of different smoothers, coarse space
configurations and parameters to solve the system matrix arising from
a purely spectral GLL discretization.

\subsection{Literature}
Spectral methods are the basic building-block of SEMs
\cite{Patera1984}, as SEMs are combinations of $h$-type Finite Element
Methods (FEM) and $p$-type spectral methods. Spectral methods can be
traced back at least to the method of selected points from Lanczos in
1956 \cite{Lanczos1956}; details of the method and its early
development can be found in \cite{GottliebOrszag1977,
  CanutoYousuffQuarteroniZang1987} and references therein. The linear
systems arising from spectral methods are dense, and inverting these
systems for a very high-order using a direct method implies using a
large memory footprint.  Direct methods offer little flexibility to a
changing geometry, accuracy and stability, for instance, during
time-stepping.  Moreover, they usually do not adapt easily to the high
levels of parallelism of today's architectures (e.g. SSE, AVX, GPU,
CPU multithreading, MPI). Iterative solvers can be a useful
alternative as long as they can achieve a rapid convergence rate and a
complexity that is bounded by the most efficient matrix-vector
multiplication technique available. For spectral methods like the ones
we use in this manuscript, using sum-factorization techniques, this
complexity bound in 2D is $\mathcal{O}\left(p^3\right)$ where
$p$ is the polynomial degree.

A variety of preconditioners have been attempted for spectral
methods. Orszag \cite{Orszag1980} used a low-order 2D finite
difference (FD) approximation.  Zang, Wong and Hussaini
\cite{ZangWongHussaini1982} introduced multigrid preconditioners for
spectral discretizations using modified Richardson iterations and
point-Jacobi smoothers. It was Brandt in his seminal paper
\cite{Brandt1977} who first suggested the idea of smoothing along
\emph{lines}, line smoothing is a special case of incomplete
factorization preconditioners spearheaded by Axelsson
\cite{Axelsson1985,Axelsson1986} who concluded that incomplete LU were
the best smoothers. Phillips \cite{Phillips1987} further evaluated
variations of Orszag's preconditioners as smoothers using relaxation
schemes in a multigrid approach.

The flexibility of the SEM attracted more research, and several
authors proposed multigrid preconditioners for SEMs around the same
time (see \cite{ZangWongHussaini1984, RonquistPatera1987,
  Ronquist1988, Heinrichs1988, Fischer1997, PavarinoWidlund1996} and
references therein), that ultimately led us to our selection of a
multigrid preconditioner for GeoFLOW. The present manuscript is part
of a ground-up development of an $hp$-multigrid for SEM, and we focus
solely on spectral methods.

We propose to draw from Phillips' work and re-evaluate the performance of
Brandt's line smoothers when combined with the multigrid $\gamma$-cycle. The
motivation for this choice is given by the fact that Golub's
cyclic-reduction \cite{BuzbeeGolubNielson70}, seen as a multigrid method, delivers a
parallelizable direct solver that can be used for line smoothers; the
use of a $\gamma$-cycle does not increase the complexity scaling of
the spectral method itself, and the overall memory footprint can be
controlled by matrix-free techniques.

In summary, we use a GMRES iterative solver and design a $p$-multigrid
\cite{ZangWongHussaini1982, ZangWongHussaini1984} $\gamma$-cycle
preconditioner with highly-vectorizable smoothers, that we refer to as
\emph{line smoothers}. Our smoothers are based on an incomplete matrix
factorization with minimal bandwidth \cite{Axelsson1985,Axelsson1986}.
The factorization delivers 1D problems that are solved by a
reformulation of cyclic reduction \cite{BuzbeeGolubNielson70} as a
direct multigrid method.

\subsection{Organization of the manuscript}
The manuscript is organized as follows: \S\ref{sec:continuoussetting}
describes the Sobolev space setting; \S\ref{sec:discretesetting}
describes the Gauss-Legendre-Lobatto discretization; and
\S\ref{sec:boundaryConsiderations} describes the treatment of boundary
conditions and the complexity of the linear system obtained.
\S\ref{sec:solver} describes the solver and preconditioners and is
divided as follows: \S\ref{sec:LinePreconditioners} describes the low
order discretizations used for smoothing; \S\ref{sec:linesolver}
describes our formulation of cyclic-reduction as multigrid to solve
the smoother linear systems; and \S\ref{sec:multigridPreconditioner}
describes the multigrid algorithm in detail. Finally,
\S\ref{sec:NumericalExperiments} provides numerical evidence of the
efficiency of the method in terms of the number of iterations required
for a fixed residual reduction.

\section{Model problem}
\subsection{Continuous Setting}\label{sec:continuoussetting}
We consider the solution of a Poisson problem in 2D: find
$u:\Omega \rightarrow \mathbb{R}$ such that
\begin{align}
  \begin{aligned}
  \mathcal{L} u := \frac{\partial^2 u}{\partial x^2} +
  \frac{\partial^2 u}{\partial y^2} = f, && &(x,y) \in \Omega \text{
    and }\\
  u = u_0, && &(x,y) \in \partial\Omega,
  \end{aligned}\label{eqn:strongFormulation}
\end{align}
where $\Omega$ is a Lipschitz domain, $\partial\Omega$ is the
intersection between $\Omega$ and its closure and $f$ is a known
function we describe later.

We introduce the Hilbert spaces $L^2(\Omega)$ and $H^1_0(\Omega)$,
where $H^1_0(\Omega)$ is the standard Sobolev space with zero boundary
trace.  They are provided with inner products $(u,v)_{L^2(\Omega)}$
and $(u,v)_{H^1_0(\Omega)} = \int_\Omega \nabla u \cdot \nabla v dx$.
The weak form of the problem reads: find $u \in H^1_0(\Omega)$ such
that
\begin{align}
  \begin{aligned}
    \mathcal{A}(u,v) := \int_\Omega \left(\frac{\partial u}{\partial x} \frac{\partial v}{\partial x} +
    \frac{\partial u}{\partial y} \frac{\partial v}{\partial y} \right) d\Omega =
    \int_{\Omega} f v d\Omega
    , && \forall v \in H^1_0(\Omega),
  \end{aligned}\label{eqn:weakFormulation}
\end{align}
where $f \in L^2(\Omega)$.  The bilinear form $\mathcal{A}(u,v)$ is
continuous and $H^1_0(\Omega)$-coercive relatively to $L^2(\Omega)$
\cite{DautrayLions1985}, i.e. there exist constants
$c_\mathcal{A},C_\mathcal{A}>0$ such that
\begin{align}
  \begin{aligned}
    \mathcal{A}(u,u) \ge c_\mathcal{A} (u,u)_{L^2(\Omega)} && \text{ and }  &&
    \mathcal{A}(u,v) \le C_\mathcal{A} (u,u)_{L^2(\Omega)}^{\frac12} (v,v)_{L^2(\Omega)}^{\frac12}.
  \end{aligned}
\end{align}
From Lax-Milgram's theorem, the variational problem admits a unique
solution in $H^1_0(\Omega)$. Non-homogeneous boundary conditions will
be considered for the discrete problem in the next section.

\subsection{Discrete setting}\label{sec:discretesetting}

We use a spectral discretization to solve system
\eqref{eqn:weakFormulation}. Our choice of spectral method is motivated by spectral
convergence, but also by the usefulness of these methods in SEMs (see
\cite{Patera1984}). Our proposed discretization and solver is used as
a smoother inside an $hp$-formulation that is ongoing work, hence we
strive for a method with a low memory footprint and high efficiency.

Several spectral method choices exist, based on Chebyshev \cite{Patera1984} and
Legendre \cite{RonquistPatera1987_2} polynomials, among others.  We
choose Gauss-Lobatto-Legendre (GLL) because of their straightforward
implementation.  A more detailed discussion on the available choices
and their characterization is available in \cite[\S
  2.4]{DevilleFischerMund2002}, see also references therein.

This section consists of two short definitions of the GLL
quadrature spanning the 2D polynomial tensor-product space that will
discretize the system \eqref{eqn:weakFormulation}. We begin with the
quadrature points and weights.
\begin{figure}
  \includegraphics{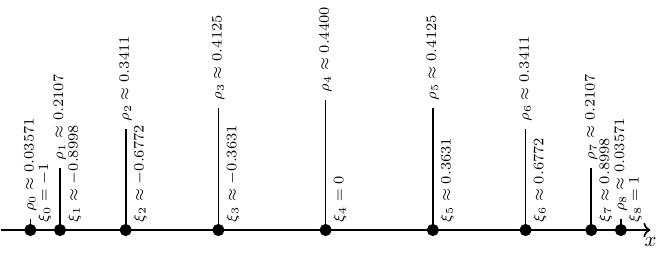}
  \caption{Quadrature points $\xi_i$ and weights $\rho_i$ for the GLL quadrature rule of order 8 (see Definition \ref{def:GLLQuadratureRule1D}).}
  \label{fig:QuadraturePoints1D}
\end{figure}
\begin{defi}[GLL quadrature rule of order $p$ (see Figure \ref{fig:QuadraturePoints1D})]
  \label{def:GLLQuadratureRule1D}
  Let $L_p(x):[-1,1]\rightarrow\mathbb{R}$ be the Legendre polynomial
  of degree $p$ \cite[p.419]{Legendre1785}. The $p+1$ GLL quadrature
  points satisfy
  \[
  \xi_0 = -1, \ \xi_p = 1, \xi_n: \ L_p'(\xi_n) = 0 \ \forall \ 0 < n < p
  \]
  and the corresponding quadrature weights are defined as
  \[
  \rho_n := \frac{2}{p(p+1)} \frac{1}{L_p^2(\xi_n)}
  .\]
  Then the GLL quadrature rule of order $p$ is defined as
  \[
  \int_{-1}^1 f(x) dx := \sum_{n=0}^p \rho_n f(\xi_n).
  \]
\end{defi}
Figure \ref{fig:QuadraturePoints1D} shows a graphical depiction of the
quadrature points and weights from Definition
\ref{def:GLLQuadratureRule1D} for polynomial degree 8, plotted on a
horizontal line and vertical lines representing the magnitude of the
quadrature weight associated with each point for the domain $[-1,1]$.
We follow by defining the interpolation basis, based on the quadrature
rule.
\begin{figure}
  \includegraphics{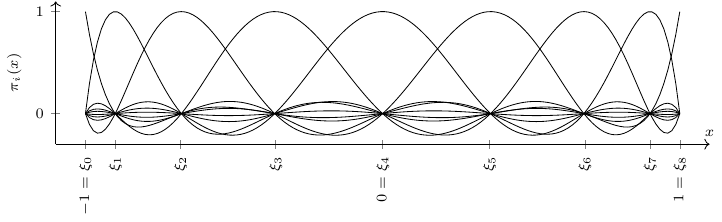}
  \includegraphics{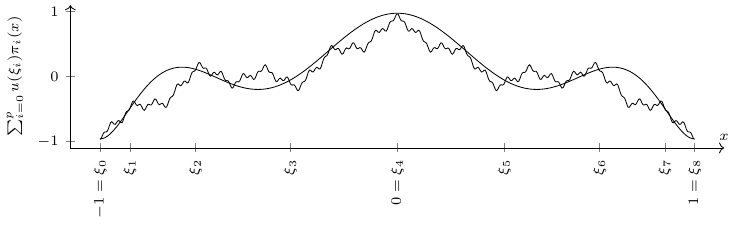}
  \caption{GLL interpolation basis (see Definition \ref{def:GLLInterpolationBasis1D}) of order 8 (top). Interpolant of a regular function $u(x)$ (bottom).}
  \label{fig:InterpolationBasis1D}
\end{figure}
\begin{defi}[GLL interpolation basis (see Figure \ref{fig:InterpolationBasis1D})]
  \label{def:GLLInterpolationBasis1D}
  The GLL interpolation basis of order $p$ is the set of Lagrange
  interpolation polynomials \cite{Lagrange1795}
  $\left\{\pi_0(x),\dots,\pi_p(x)\right\}$ with roots
  $\left\{\xi_0,\dots,\xi_p\right\}$.  Let
  $u(x):\mathbb{R}\rightarrow\mathbb{R}$ be a regular function. Its
  interpolant is written as
  \[u(x) \approx \sum_{i=0}^p u\left(\xi_i\right) \pi_i(x),\]
  see \cite[Eq. (2.4.3)]{DevilleFischerMund2002}.
\end{defi}
Figure \ref{fig:InterpolationBasis1D} shows the polynomial basis from
Definition \ref{def:GLLInterpolationBasis1D} for polynomial degree 8
and an example of the polynomial interpolant obtained for an arbitrary
function $u(x)$. We continue with the span that allows us to get a
basis for a 2D space.
\begin{figure}
  \begin{subfigure}{0.24\textwidth}
  \includegraphics{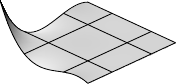}
    \caption{$v_{1}(x,y)$}
  \end{subfigure}
  \begin{subfigure}{0.24\textwidth}
  \includegraphics{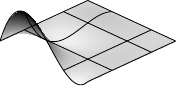}
    \caption{$v_{2}(x,y)$}
  \end{subfigure}
  \begin{subfigure}{0.24\textwidth}
  \includegraphics{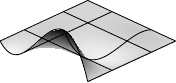}
    \caption{$v_{3}(x,y)$}
  \end{subfigure}
  \begin{subfigure}{0.24\textwidth}
  \includegraphics{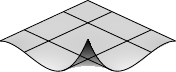}
    \caption{$v_{4}(x,y)$}
  \end{subfigure}

  \begin{subfigure}{0.24\textwidth}
  \includegraphics{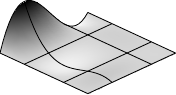}
    \caption{$v_{5}(x,y)$}
  \end{subfigure}
  \begin{subfigure}{0.24\textwidth}
  \includegraphics{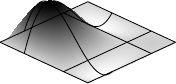}
    \caption{$v_{6}(x,y)$}
  \end{subfigure}
  \begin{subfigure}{0.24\textwidth}
  \includegraphics{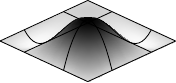}
    \caption{$v_{7}(x,y)$}
  \end{subfigure}
  \begin{subfigure}{0.24\textwidth}
  \includegraphics{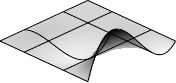}
    \caption{$v_{8}(x,y)$}
  \end{subfigure}

  \begin{subfigure}{0.24\textwidth}
  \includegraphics{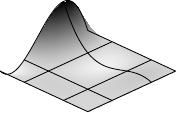}
    \caption{$v_{9}(x,y)$}
  \end{subfigure}
  \begin{subfigure}{0.24\textwidth}
  \includegraphics{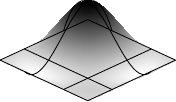}
    \caption{$v_{10}(x,y)$}
  \end{subfigure}
  \begin{subfigure}{0.24\textwidth}
  \includegraphics{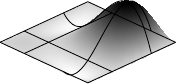}
    \caption{$v_{11}(x,y)$}
  \end{subfigure}
  \begin{subfigure}{0.24\textwidth}
  \includegraphics{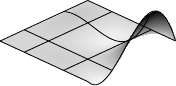}
    \caption{$v_{12}(x,y)$}
  \end{subfigure}
 
  \begin{subfigure}{0.24\textwidth}
  \includegraphics{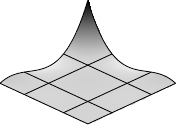}
    \caption{$v_{13}(x,y)$}
  \end{subfigure}
  \begin{subfigure}{0.24\textwidth}
  \includegraphics{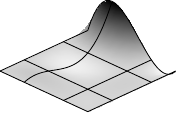}
    \caption{$v_{14}(x,y)$}
  \end{subfigure}
  \begin{subfigure}{0.24\textwidth}
  \includegraphics{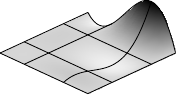}
    \caption{$v_{15}(x,y)$}
  \end{subfigure}
  \begin{subfigure}{0.24\textwidth}
  \includegraphics{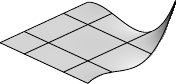}
    \caption{$v_{16}(x,y)$}
  \end{subfigure}
  \begin{subfigure}{0.49\textwidth}
  \includegraphics{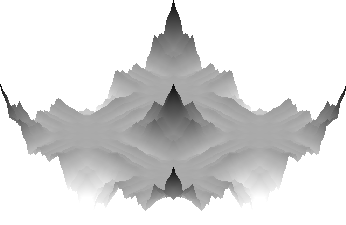}
    \caption{$u(x,y)$}
  \end{subfigure}
  \begin{subfigure}{0.49\textwidth}
  \includegraphics{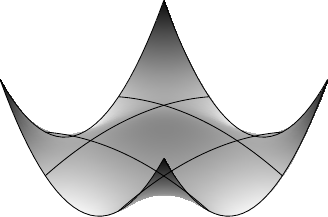}
    \caption{$\sum_{k=0}^{p} \sum_{j=0}^{p} u(\mathbf{x}_{kl}) v_{\left((p+1)k+(j+1)\right)}(x,y)$}
  \end{subfigure}
  \caption{GLL tensor product interpolation basis of order 3 on a square domain (Figures A to P).
    Order 3 interpolant of a regular function $u(x,y)$ defined in a square domain (Figures Q and R).}
\end{figure}

We have all the elements needed to generate the space of 2D basis
functions needed to discretize our problem. Let
\[\mathbb{P}_p\left([-1,1,]^2\right) := \text{span}\left\{
\left\{\pi_0(x),\dots,\pi_p(x)\right\} \otimes
\left\{\pi_0(y),\dots,\pi_p(y)\right\} \right\}\] be the space of
tensor product GLL interpolant polynomials of degree $p$ in $x$ and
$y$ on $[-1,1]^2$, and let $\mathbb{T}$ be a domain described by a
mapping $\Phi:[-1,1]^2 \rightarrow \mathbb{T}$. Assuming that the
Jacobian of $\Phi$ and its inverse is uniformly bounded, define
the mapped space $\mathbb{P}_p(\mathbb{T})$ as the pull-back of
functions under $\Phi$. We define the space $V_p$ as
\begin{equation}
  V_p := \left\{v \in L^2(\mathbb{T}) \big| v \in
  \mathbb{P}_p(\mathbb{T}) \right\}
  \label{eqn:DiscreteSpace}
\end{equation}
where the GLL quadrature rule is used to calculate integrals and inner
products, such that for $u \in V_p$
\begin{align*}
  \int_\mathbb{T} u~d\mathbf{x} := \sum_{k=0}^{p} \sum_{l=0}^{p} \phi(\mathbf{x}_{kl}) u(\mathbf{x}_{kl}),
\end{align*}
where $\mathbf{x}_{kl} = \Phi((\xi_k,\xi_l)), \forall k,l \in
\{0,\dots,p\}$ are the 2D quadrature points and $\phi(\mathbf{x}_{ij})$
are the 2D quadrature weights that include geometric terms to account
for deformations. The methodology to include deformations exceeds the
scope of this manuscript and can be found explained in detail in
\cite[p. 180]{DevilleFischerMund2002} and \cite[p. 26]{Ronquist1988}.

With these definitions, the entries of the mass matrix $M$ and the
Laplace stiffness matrix $A$ are
\begin{align}
  \left\{
  \begin{aligned}
    M_p: \left(M_{ij}\right) =& \int_\mathbb{T} v_i v_j d\mathbf{x} = \sum_{k=0}^p \sum_{l=0}^p \phi(\mathbf{x}_{kl}) v_i(\mathbf{x}_{kl}) v_j(\mathbf{x}_{ij}) \\
    A_p: \left(A_{ij}\right) =& \int_\mathbb{T} \left(\frac{\partial v_i}{\partial x} \frac{\partial v_j}{\partial x} + \frac{\partial v_i}{\partial y} \frac{\partial v_j}{\partial y}\right) d\mathbf{x} \\
    =& \sum_{k=0}^p \sum_{l=0}^p \phi(\mathbf{x}_{kl}) \left(\frac{\partial v_i}{\partial x} (\mathbf{x}_{kl}) \frac{\partial v_j}{\partial x}(\mathbf{x}_{kl}) + \frac{\partial v_i}{\partial y}(\mathbf{x}_{kl}) \frac{\partial v_j}{\partial y}(\mathbf{x}_{kl}) \right)
  \end{aligned}\right.
  \label{eqn:discreteBilinear}
\end{align}
where $\forall v_i \in V_p$, such that the discretized linear system can be written as
\begin{equation}
  A_p \mathbf{u} = M_p \mathbf{f},
  \label{eqn:discreteSystem}
\end{equation}
where $\mathbf{u}$ and $\mathbf{f}$ are the coefficient vectors of the
representation of $u$ and $f$, in terms of the chosen
basis which by its definition are also the values of the functions at
the quadrature points.

The matrix $A_p$ has convenient properties in terms of the cost of
applying it to a vector that make spectral methods competitive with
other, lower order methods \cite{Roth1934}. We refer the reader to
\cite[\S 4]{DevilleFischerMund2002} for details on how the tensor
product structure is kept in the presence of deformations.

\subsection{Boundary considerations}\label{sec:boundaryConsiderations}
The enforcement of boundary conditions follows the well-established approach
of pre- and post-processing of the discrete right hand side. We first outline this
method in the continuous setting before detailing its discrete
implementation.

Let $u$ be the solution of problem \ref{eqn:strongFormulation} and let
$v$ be an arbitrary function such that
\[
v\big|_{\partial\Omega} := u\big|_{\partial\Omega}.
\]
Defining $w:=u-v$, we have
\begin{align*}
  w\big|_{\partial\Omega} =& (u-v)\big|_{\partial\Omega} = 0 \text{ and }\\
  \mathcal{L} w =& \mathcal{L}(u-v) = \mathcal{L} u - \mathcal{L} v = f - \mathcal{L} v.
\end{align*}
Choosing $v=0$ in $\Omega$ we have $\mathcal{L}v=0$, thus we
can first solve the problem: Find $w$ such that
\begin{align*}
  \left\{
  \begin{aligned}
    \mathcal{L} w =& f \text{ in } \Omega\\
    w\big|_{\partial\Omega} =& 0
  \end{aligned}\right. ,
\end{align*}
and then obtain the solution $u$ as $u = w + v$.

In discrete form, we have
\begin{align*}
  \begin{aligned}
    A_p \mathbf{u} =& M_p \mathbf{f} \\
    A_p \left(\mathbf{w} + \mathbf{v}\right) =& M_p \mathbf{f} \\
    A_p \mathbf{w} =& M_p \mathbf{f} - A_p \mathbf{v} := \mathbf{f}_0,
  \end{aligned}
\end{align*}
thus, we solve the problem $A_p \mathbf{w} = \mathbf{f}_0$ and we
obtain the solution for non-homogeneous boundary conditions as
$\mathbf{u} =\mathbf{w} + \mathbf{v}$.

Given the tensor product structure of the matrix $A_p$, the cost of
applying it to a vector is $\mathcal{O}\left(p^3\right)$ using the
fast application of tensor products \cite{Roth1934}. For structured
domains the system could be solved by a fast diagonalization method
with complexity $\mathcal{O}\left(p^3\right)$, but in our work we
require a low memory footprint and solutions on deformed geometries
\cite[\S 4.4]{DevilleFischerMund2002}; therefore, we propose an
iterative solution of the system as in \cite[\S
  4.5.4]{DevilleFischerMund2002}, using a dedicated
$\mathcal{O}\left(p^3\right)$ multigrid preconditioner that will be
described in the next section.

\section{Solver}\label{sec:solver}
The system matrix we will invert, $A_p$, is symmetric but no assumption
will be made on the symmetry of the system matrix for the
applicability of the solvers and preconditioners presented in this
manuscript. For instance, a simple change to Chebyshev polynomials
instead of Legendre polynomials would deliver a nonsymmetric matrix.
Thus, we have chosen GMRES as a solver, instead of the classical
choice of CG for symmetric problems.

For scalability purposes, since the application of our system matrix
is $\mathcal{O}\left(p^3\right)$, we would like to maintain the
complexity; therefore, we need a preconditioner that would allow the
iteration count to be bounded as the problem becomes larger, while at the
same time costing no more than $\mathcal{O}\left(p^3\right)$. For this
purpose, we use a multigrid preconditioner with direct multigrid line
smoothers that we will describe in the next sections and which
constitutes the heart of our method.  We will describe the elements
that make up the multigrid procedure in separate subsections.

\subsection{Line preconditioners}\label{sec:LinePreconditioners}
We choose line preconditioners, for which well-known theory is
available (see \cite{Axelsson1985,Axelsson1986} and references
therein). This section defines the structure of our preconditioners
that will be used afterwards to design multilevel smoothers. The
definition of these preconditioners is algebraic in nature and has
been spearheaded by Axelsson. We reproduce
\cite[Def. 2.1]{Axelsson1985} below, which in our manuscript we use
only with half-bandwidth 1.
\begin{defi}[Axelsson incomplete bandwidth preconditioner]
  Let $H$ be a square matrix and let $p\ge 0$ be an integer, then
  $[H]^{(p)}$ denotes the matrix with entries equal to those within
  the band position of $H$ with half-bandwidth $p$ and zero outside,
  i.e.,
  \[[H]_{i,j}^{(p)}=\begin{cases} H_{i,j}, &|i-j| \le p, \\
  0, &\text{otherwise.} \end{cases}\]
  \label{def:halfbandwidth}
\end{defi}
\noindent The main interest of these preconditioners is that they are strongly
vectorizable and applicable straightforwardly in a matrix-free fashion
while using SSE and AVX instructions in modern microprocessors. A
subset of the elements of $A_p$ are chosen and kept, and the resulting
system is inverted using a fast method.

Literature shows other approaches involving overlaying a first order
finite-difference or finite-element method on the quadrature point
mesh \cite{Orszag1980,DevilleMund1985,PartnerRothman1995,Olson2007},
so we also use the line preconditioner approach over a matrix
$\tilde{A}_p$ obtained as the discretization with bilinear finite
elements of the system \eqref{eqn:weakFormulation} on the GLL quadrature
points.

We will define the preconditioner \emph{systems} algebraically first, 
before detailing the solver used to invert them.

\begin{figure}
  \centering
  \begin{subfigure}{0.42\textwidth}
  \resizebox{\textwidth}{!}{ % Adjust the size to the text width
  \includegraphics{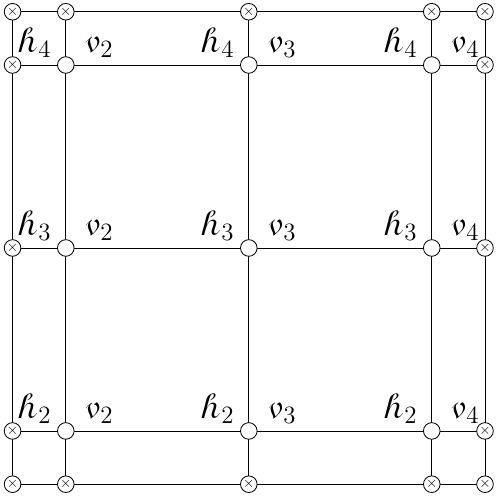}
  }
  \caption{Mesh}
  \end{subfigure}
  \begin{subfigure}{0.45\textwidth}
    \begin{align*}
      A_4 =& \left(\vcenter{\hbox{\resizebox{\textwidth}{!}{
            \includegraphics{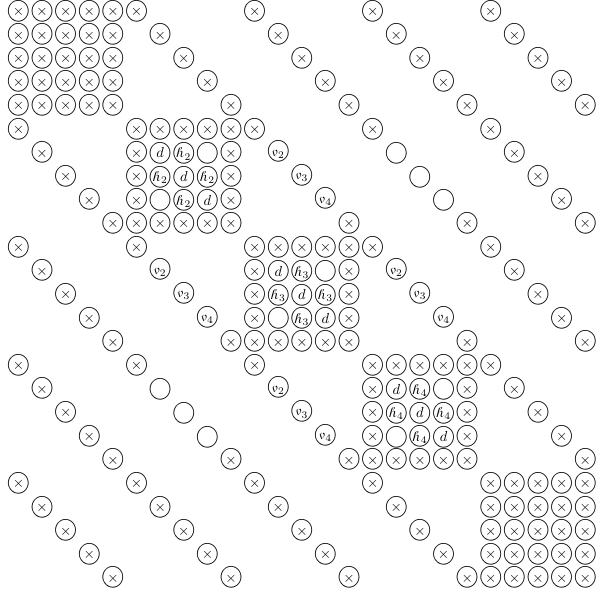}}}}\right)
    \end{align*}
    \caption{GLL Matrix sparsity}
  \end{subfigure}
  \begin{subfigure}{0.45\textwidth}
    \begin{align*}
      \tilde{A}_4 =& \left(\vcenter{\hbox{\resizebox{\textwidth}{!}{
            \includegraphics{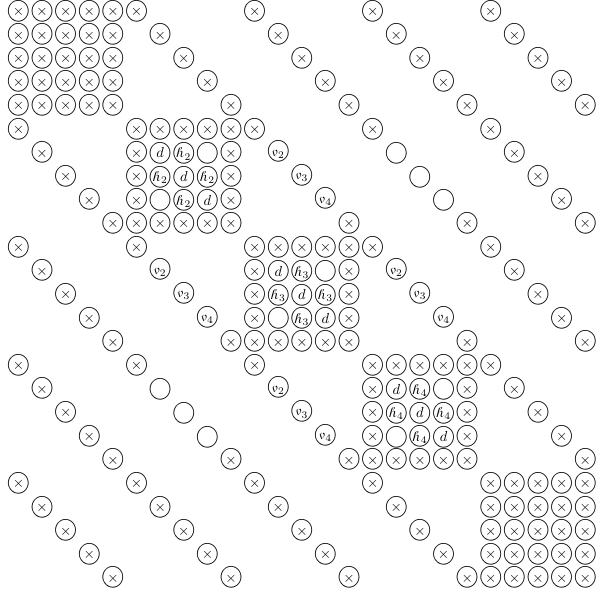}}}}\right)
    \end{align*}
    \caption{FEM Matrix sparsity}
  \end{subfigure}
  \caption{
    Elements associated to each line of the line smoother
    for $p=4$. \\
    $\rightarrow$ Elements represented by a circle \textcircled{}
    (marked or unmarked) are the non-zero elements of the GLL or FEM
    discretized Laplacian, using lexicographic node ordering
    left-right, bottom-top,\\
    $\rightarrow$ circles marked by \scalebox{0.5}{$\times$} are not
    needed due to boundary conditions,\\
    $\rightarrow$ circles marked by $d$ are always part of the line
    systems they share a row or column with, whether it's $\mathscr{h}_i$ or
    $\mathscr{v}_i$, for all $i=1,\dots,p-2$,\\
    $\rightarrow$ circles marked by $\mathscr{h}_i$, together with the
    diagonal for each line, conform the system for each horizontal
    line $i=1,\dots,p-2$, and \\
    $\rightarrow$ circles marked by $\mathscr{v}_i$, together
    with the diagonal for each line, conform the system for each
    vertical line $i=1,\dots,p-2$.
  }
  \label{fig:LinePreconditioners}
\end{figure}

Figure \ref{fig:LinePreconditioners}
shows a GLL mesh (left) for polynomial degree $4$ with nodes marked
following horizontal lines $\mathscr{h}_i$ and vertical lines
$\mathscr{v}_i$, and the sparsity pattern associated with the GLL
discretization (right), following a lexicographic left-right, then
bottom-top, node ordering.

It is clear that if we apply Def. \ref{def:halfbandwidth} with
half-bandwidth equal to $1$, we obtain the systems noted by
$\mathscr{h}_i$ that can be solved in parallel since they form a
block-diagonal matrix. In the same fashion, should we choose a
lexicographic ordering that first runs bottom-top and then left-right,
the systems $\mathscr{v}_i$ will also be block-diagonal and can be
solved in parallel.

We will call these systems \emph{line systems} and by solving them in
parallel we obtain \emph{line preconditioners} that boil down to
block-Jacobi preconditioners for the two lexicographic orderings
considered.
\begin{defi}[Line preconditioners]
  Let $R_{\mathscr{h}_i}:\mathbb{R}^{(p^2)}\rightarrow\mathbb{R}^{(p-1)}$
  (respectively $R_{\mathscr{v}_i}$) be the operator that extracts the
  degrees of freedom of each $i=2,\cdots,p$ horizontal (respectively
  vertical) line in the GLL discretization, and sets all other elements
  to zero. The {\bf horizontal line preconditioner} is defined as
  \[B_{\mathscr{h}}^{-1} = \sum_{i=2}^{p} P_{\mathscr{h}_i} \left(R_{\mathscr{h}_i} A_p P_{\mathscr{h}_i}\right)^{-1} R_{\mathscr{h}_i}, \]
  where $P_{\mathscr{h}_i}=R_{\mathscr{h}_i}^\intercal$, and
  the {\bf vertical line preconditioner} is
  \[B_{\mathscr{v}}^{-1} = \sum_{i=2}^{p} P_{\mathscr{v}_i} \left(R_{\mathscr{v}_i} A_p P_{\mathscr{v}_i}\right)^{-1} R_{\mathscr{v}_i}, \]
  where $P_{\mathscr{v}_i}=R_{\mathscr{v}_i}^\intercal$.

  The {\bf FEM line preconditioners} $\tilde{B}_{\mathscr{h}}^{-1}$
  and $\tilde{B}_{\mathscr{v}}^{-1}$ are defined analogously, simply
  replacing $A_p$ by the matrix $\tilde{A}_p$, resulting from the
  bilinear finite element discretization of
  \eqref{eqn:weakFormulation} on a mesh consisting of the GLL
  quadrature nodes.
\end{defi}
Condition number estimates for the application of this preconditioner
can be found in \cite{Demko1984} for arbitrary bandwidths.
It is clear, however, that for our half-bandwidth choice, the system
can be inverted line-by-line by solving tridiagonal systems. For
this task, we choose a specific formulation of the well known
cyclic-reduction solver \cite{BuzbeeGolubNielson70}.

\subsection{Line solver}\label{sec:linesolver}
We will focus on the local solver used to solve each system
$\left(R_{\mathscr{h}_i} A_p P_{\mathscr{h}_i}\right)^{-1}$ and
$\left(R_{\mathscr{v}_i} A_p P_{\mathscr{v}_i}\right)^{-1}$, for
$i=2,\dots,p$. To the best of our knowledge, even though the
literature contains many implementations of block-cyclic reduction, a
formal proof of it formulated as a special case of multigrid has not
been made available.

We describe a generic, parallel solver for block-tridiagonal systems
that has its roots in the well-known cyclic reduction algorithm of
Buzbee, Golub and Nielson \cite{BuzbeeGolubNielson70}. We reformulate
cyclic reduction as a multigrid V-cycle without post-smoothing in
order to employ the typical data structures of multigrid for its
implementation. Our reformulation underlines the places where
parallelism can be applied and brings together some of the findings
from Gander, Kwok and Zhang \cite{GanderKwokZhang2018} as well.

The understanding of cyclic-reduction as a two-level preconditioned
solver explains the spurious appearance of exact solvers in Local
Fourier Analyses such as \cite{GanderLucero22_2,LuceroLorcaGander2024}. It
shows that when optimizing 2D and 3D two-level solvers by using their
1D version (see \cite{Hemker2003,Hemker2004} and references therein),
if we are allowed to choose different restriction and prolongation
operators, the optimum boils down to the restriction and prolongation
operators that we describe herein, but we will make the complete link
with Local Fourier Analysis of a diverse set of discretizations in an
upcoming separate manuscript.

Proposition \ref{prop:blockinverse} shows that the well-known formula
to invert a $2\times 2$-block matrix (a simple reformulation of
\cite[Corollary 2.8.9]{Bernstein2009}), can be rearranged to a form
that exposes the sparsity of the elements of the formula and
underlines the re-usage of the data structures involved.
\begin{prop}
  Let $M \in \mathbb{R}^{n \times n}$, $n \in \mathbb{Z}^+$ have the
  block structure
  \begin{equation*}
    M = \left(
    \begin{array}{cc}
      A & B \\
      C & D
    \end{array}
    \right),
  \end{equation*}
  assuming $D-C A^{-1} B$ is invertible, the inverse of $M$ can be
  expressed as
  \begin{align*}
    \scalemath{0.8}{
      M^{-1} =
      \left(
      \begin{array}{cc}
        A^{-1} & 0 \\
        0 & 0
      \end{array}
      \right)
      +
      \left(
      \begin{array}{c}
        -A^{-1} B \\
        I
      \end{array}
      \right)
      (D - C A^{-1} B)^{-1}
      \left(
      \begin{array}{cc}
        -\left(A^{-1}B\right)^\intercal & I
      \end{array}
      \right)
      \left(
      I
      -
      M
      \left(
      \begin{array}{cc}
        A^{-1} & 0 \\
        0 & 0
      \end{array}
      \right)
      \right)  },
  \end{align*}
  where $I$ is the identity matrix of appropriate size and the left
  block of $\left(\begin{array}{cc}-\left(A^{-1}B\right)^\intercal &
    I\end{array}\right)$ is arbitrary.
  \label{prop:blockinverse}
\end{prop}
\begin{proof}
  The proof is straightforward, starting with a simple reformulation of a
  known formula that can be found in \cite[Corollary
    2.8.9]{Bernstein2009}
  \begin{align*}
    M^{-1}
    =& \begin{pmatrix}
      A^{-1} + A^{-1} B(D - C A^{-1} B)^{-1} C A^{-1} & -A^{-1} B(D - C A^{-1} B)^{-1} \\
      -(D - C A^{-1} B)^{-1} C A^{-1} & (D - C A^{-1} B)^{-1}
    \end{pmatrix} \\
    =& \begin{pmatrix}
      A^{-1} & 0 \\
      0 & 0
    \end{pmatrix}
    +
    \begin{pmatrix}
      A^{-1} B(D - C A^{-1} B)^{-1} C A^{-1} & -A^{-1} B(D - C A^{-1} B)^{-1} \\
      -(D - C A^{-1} B)^{-1} C A^{-1} & (D - C A^{-1} B)^{-1}
    \end{pmatrix} \\
    =& \begin{pmatrix}
      A^{-1} & 0 \\
      0 & 0
    \end{pmatrix}
    +
    \begin{pmatrix}
      -A^{-1} B \\
      I
    \end{pmatrix}
    (D - C A^{-1} B)^{-1}
    \begin{pmatrix}
      -C A^{-1} & I
    \end{pmatrix} \\
    =& \begin{pmatrix}
      A^{-1} & 0 \\
      0 & 0
    \end{pmatrix}
    +
    \begin{pmatrix}
      -A^{-1} B \\
      I
    \end{pmatrix}
    (D - C A^{-1} B)^{-1}
    \begin{pmatrix}
      -\left(A^{-1}B\right)^\intercal & I
    \end{pmatrix}
    \begin{pmatrix}
      0 & 0 \\
      -C A^{-1} & I
    \end{pmatrix},
  \end{align*}
  where $I$ is the identity matrix of appropriate size and the left
  block of $\left(\begin{array}{cc}-\left(A^{-1}B\right)^\intercal &
    I\end{array}\right)$ is arbitrary since it gets post-multiplied by
    a zero block.  The result is achieved with some final
    manipulation:
  \begin{align*}
    \begin{pmatrix} 0 & 0 \\ -C A^{-1} & I \end{pmatrix}
    = & I - \begin{pmatrix} I & 0 \\ C A^{-1} & 0 \end{pmatrix} \\
    = & I - M \begin{pmatrix} A^{-1} & 0 \\ 0 & 0 \end{pmatrix}.
  \end{align*}
\end{proof}
We reinterpret Proposition \ref{prop:blockinverse} as a two-level
preconditioned Richardson solver. For clarity, we describe the
preconditioned solver in Algorithm
\ref{alg:twolevelpreconditionernosmoothing}.
\begin{algorithm}
  \caption{Two-level multigrid preconditioner $(S,M_0)$ without post-smoothing. Define
  the action of the operator $M^{-1}$ on a vector $\boldsymbol{g}$ as:}
  \begin{algorithmic}[1]
    \State compute $\boldsymbol{x}:= S^{-1} \boldsymbol{g}$,
    \State compute $\boldsymbol{y}:= \boldsymbol{x} + P M_0^{-1} R
    (\boldsymbol{g} - M \boldsymbol{x})$,
    \State obtain $M^{-1}\boldsymbol{g} = \boldsymbol{y}$.
  \end{algorithmic}
  \label{alg:twolevelpreconditionernosmoothing}
\end{algorithm}

Lemma \ref{lem:inversionrichardsonequivalence} draws the fundamental
link between the $2\times 2$-block inversion formula and a two-level
preconditioned Richardson solver without post-smoothing.
\begin{lem}[Block inversion equivalence to two-level preconditioned Richardson without post-smoothing]
  Let
  \begin{multline*}
    S^{-1} =
    \left(
    \begin{array}{cc}
      A^{-1} & 0 \\
      0 & 0
    \end{array}
    \right), 
    P = 
    \left(
    \begin{array}{c}
      -A^{-1} B \\
      I
    \end{array}
    \right), 
    R = P^\intercal
    \text{ and }
      M_0 = D - C A^{-1} B.
  \end{multline*}
  With these definitions, the application of Proposition
  \ref{prop:blockinverse} is equivalent to one iteration of
  Richardson, preconditioned with Algorithm
  \ref{alg:twolevelpreconditionernosmoothing}.
  \label{lem:inversionrichardsonequivalence}
\end{lem}
\begin{proof}
  Consider a single Richardson iteration to solve a system of the form
  $M \boldsymbol{u} = \boldsymbol{f}$. With a zero initial guess, the
  residual is $\boldsymbol{g} = \boldsymbol{f}$. Apply Algorithm
  \ref{alg:twolevelpreconditionernosmoothing}, the first step is
  \[ \boldsymbol{x} = S^{-1}\boldsymbol{g}, \]
  the second step is,
  \[ \boldsymbol{y} = S^{-1}\boldsymbol{g} + P M_0^{-1} R  (\boldsymbol{g} - M S^{-1}\boldsymbol{g})
  = \left(S^{-1} + P M_0^{-1} R  (I - M S^{-1})\right) \boldsymbol{g}, \]
  and finally we obtain
  \[ M^{-1} = S^{-1} + P M_0^{-1} R  (I - M S^{-1}), \]
  where we identify the result in Proposition \ref{prop:blockinverse}
  with the given definitions of $S^{-1}$, $P$, $R$ and $M_0$.
\end{proof}

We can therefore draw a link between a preconditioned Richardson
iteration and a direct solver that is suggested in
\cite{GanderKwokZhang2018} by Corollary
\ref{cor:inversionequivalence}.

This concludes the treatment of $2\times 2$ block-matrix solvers.  We
continue with the reordering necessary to make block-tridiagonal
matrices $2\times 2$ block-matrices. To that end, we define the
notation for block tridiagonal matrices in Definition
\ref{def:blocktridiagonalmatrix}, and the general structure of the
reordering matrix $Q$ in Definition \ref{def:oddevenreordermatrix}.
\begin{defi}[Block-tridiagonal matrix]
Let $M$ be a block-tridiagonal matrix with the following notation for
the blocks:
\begin{equation}
  M=
\left(
\begin{array}{ccccccccc}
    A_{1} &  U_{1} &     0 &     0 & \dots & 0 & 0 & 0 & 0\\
    L_{2} &  A_{2} & U_{2} &     0 & \dots & 0 & 0 & 0 & 0\\
        0 &  L_{3} & A_{3} & U_{3} & \dots & 0 & 0 & 0 & 0\\
        0 &      0 & L_{4} & A_{4} & \dots & 0 & 0 & 0 & 0\\
    \dots &  \dots & \dots & \dots & \dots & \dots & \dots & \dots & \dots\\
        0 &      0 &     0 &     0 & \dots & A_{n-3} & U_{n-3} &       0 & 0 \\
        0 &      0 &     0 &     0 & \dots & L_{n-2} & A_{n-2} & U_{n-2} & 0 \\
        0 &      0 &     0 &     0 & \dots & 0 & L_{n-1} & A_{n-1} & U_{n-1} \\
        0 &      0 &     0 &     0 & \dots & 0 &       0 & L_{n} & A_{n}
\end{array}
\right).
\end{equation}
\label{def:blocktridiagonalmatrix}
\end{defi}

\begin{defi}[Odd-even column reordering matrix.]
  Let $Q$ be the matrix that groups all the columns with an odd index
  first and all the columns with an even index afterwards, i.e.
  \begin{equation*}
    Q = \begin{pmatrix} I_{\text{odd}} & I_{\text{even}} \end{pmatrix},
  \end{equation*}
  where $I_{\text{odd}}$ represents the odd columns of the identity
  matrix, and $I_{\text{even}}$, the even columns.
  \label{def:oddevenreordermatrix}
\end{defi}

Proposition \ref{prop:linecoarsespace} shows that $Q$ can indeed
transform block-tridiagonal matrices into $2\times 2$ block matrices
with a convenient structure that simplifies the inversions required to
use a $2\times 2$ inversion formula. This is a consequence of $M$ having
block-wise ``Property A'' (see \cite{Young1954,Young1972}).
\begin{prop}
  Let $M \in \mathbb{R}^{n \times n}$, $n \in \mathbb{Z}^+$ be a
  block-tridiagonal matrix as in Definition
  \ref{def:blocktridiagonalmatrix} and $Q$ be an odd-even reordering
  matrix as in Definition \ref{def:oddevenreordermatrix}.  The matrix
  $Q^\intercal M Q$ can be divided into 4 blocks as
  \begin{equation*}
     Q^\intercal M Q = \left(
    \begin{array}{cc}
      A & B \\
      C & D
    \end{array}
    \right),
  \end{equation*}
  where the matrix $A^{-1} B$ is banded with the stencil
  \begin{align*}
    \left[\left(A_{i}^{-1} L_{i}\right)\quad\left(A_{i}^{-1} U_{i}\right)\quad\left(0\right)\right], \quad i=1,3,5,\dots, 
  \end{align*}
  and the matrix $D - C A^{-1} B$ is block-tridiagonal with the stencil
  \begin{multline*}
    \left[\left(-L_i A_{i-1}^{-1} L_{i-1}\right)\quad\left(-L_i A_{i-1}^{-1} L_{i-1} + A_i - U_i A_{i+1}^{-1} L_{i+1}\right)\quad\left(-U_i A_{i+1}^{-1} U_{i+1}\right)\right], \\ \quad i=2,4,6,\dots.
  \end{multline*}
  \label{prop:linecoarsespace}
\end{prop}
\begin{proof}
  The calculation is straightforward,
  \begin{align*}
    Q^\intercal M Q=&
    \begin{aligned}
    \left(
    \begin{array}{cc}
      \left[
        \begin{array}{ccccc}
          A_{1} &      0 &     0 &     0 & \dots \\
          0 &  A_{3} &     0 &     0 & \dots \\
          0 &      0 & A_{5} &     0 & \dots \\
          0 &      0 &     0 & A_{7} & \dots \\
          \dots &  \dots & \dots & \dots & \dots \\
        \end{array}
        \right] &
      \left[
        \begin{array}{ccccc}
          U_{1} &     0 &     0 &     0 & \dots\\
          L_{3} & U_{3} &     0 &     0 & \dots\\
          0 & L_{5} & U_{5} &     0 & \dots\\
          0 &     0 & L_{7} & U_{7} & \dots\\
          \dots & \dots & \dots & \dots & \dots\\
        \end{array}
        \right] \\
      \left[
        \begin{array}{ccccc}
          L_{2} & U_{2} &     0 &     0 & \dots\\
          0 & L_{4} & U_{4} &     0 & \dots\\
          0 &     0 & L_{6} & U_{6} & \dots\\
          0 &     0 &     0 & L_{8} & \dots\\
          \dots & \dots & \dots & \dots & \dots\\
        \end{array}
        \right] &
      \left[
        \begin{array}{ccccc}
          A_{2} &      0 &     0 &     0 & \dots \\
          0 &  A_{4} &     0 &     0 & \dots \\
          0 &      0 & A_{6} &     0 & \dots \\
          0 &      0 &     0 & A_{8} & \dots \\
          \dots &  \dots & \dots & \dots & \dots \\
        \end{array}
        \right]
    \end{array}
    \right) 
  \end{aligned}, \\
    A^{-1}B =& 
    \begin{aligned}
      \left[
        \begin{array}{ccccc}
          A_{1}^{-1} U_{1} &     0 &     0 &     0 & \dots\\
          A_{3}^{-1}L_{3} & A_{3}^{-1}U_{3} &     0 &     0 & \dots\\
          0 & A_{5}^{-1}L_{5} & A_{5}^{-1} U_{5} &     0 & \dots\\
          0 &     0 & A_{7}^{-1}L_{7} & A_{7}^{-1} U_{7} & \dots\\
          \dots & \dots & \dots & \dots & \dots\\
        \end{array}
        \right],
    \end{aligned} \\
    C A^{-1} B =&
    \scalemath{0.7}{
    \begin{aligned}
      \left[
        \begin{array}{ccccc}
          L_2 A_{1}^{-1} U_{1} + U_2 A_{3}^{-1}L_{3} & U_2 A_{3}^{-1} U_{3} &     0 &     0 & \dots\\
          L_4 A_{3}^{-1} L_{3} & L_4 A_{3}^{-1}U_{3} + U_4 A_{5}^{-1} L_{5} &  U_4 A_{5}^{-1}U_{5} &     0 & \dots\\
          0 & L_6 A_{5}^{-1} L_{5} & L_6 A_{5}^{-1} U_{5} + U_6 A_{7}^{-1} L_{7} & U_6 A_{7}^{-1}U_{7} & \dots\\
          0 &     0 & L_8 A_{7}^{-1}L_{7} & L_8 A_{7}^{-1} U_{7} + U_8 A_{9}^{-1} L_{9} & \dots\\
          \dots & \dots & \dots & \dots & \dots\\
        \end{array}
        \right]
    \end{aligned}
    },
  \end{align*}
  and the result is achieved.
\end{proof}

Lemma \ref{lem:cyclicequivalence} links cyclic-reduction to
multigrid \cite[Lemma 1]{GanderKwokZhang2018}, and the reordering matrix $Q$.

We have all the tools needed to write the general structure of choice
of prolongation and restrictions that gives a direct two-level
preconditioned solver without post-smoothing, Theorem
\ref{theo:twoleveldirect}.
\begin{theorem}[Two-level direct solver]
  Let $M$ be a block-tridiagonal matrix. A Richardson iteration
  preconditioned with Algorithm
  \ref{alg:twolevelpreconditionernosmoothing} with the 
  definitions
  \begin{align}
    S^{-1} =&
    \begin{aligned}
      \left[
        \begin{array}{ccccccccc}
          A_{1}^{-1} &      0 & 0 &     0 & 0 &     0 & 0 & \dots \\
          0 & 0 &  0 & 0 &     0 & 0 &     0 & \dots \\
          0 & 0 &  A_{3}^{-1} &     0 & 0 &     0 & 0 & \dots \\
          0 & 0 &  0 & 0 &     0 & 0 &     0 & \dots \\
          0 & 0 &      0 & 0 & A_{5}^{-1} &     0 & 0 & \dots \\
          0 & 0 &  0 & 0 &     0 & 0 &     0 & \dots \\
          0 & 0 &      0 & 0 &     0 & 0 & A_{7}^{-1} & \dots \\
          \dots &  \dots & \dots & \dots & \dots & \dots & \dots & \dots\\
        \end{array}
        \right]
      \label{eqn:lineSolverSmoother}
    \end{aligned} \\
    P =&
    \begin{aligned}
      \left[
        \begin{array}{ccccc}
          -A_{1}^{-1} U_{1} &     0 &     0 &     0 & \dots\\
          I &     0 &     0 &     0 & \dots\\
          -A_{3}^{-1}L_{3} & -A_{3}^{-1}U_{3} &     0 &     0 & \dots\\
          0 &     I &     0 &     0 & \dots\\
          0 & -A_{5}^{-1}L_{5} & -A_{5}^{-1} U_{5} &     0 & \dots\\
          0 &     0 &     I &     0 & \dots\\
          0 &     0 & -A_{7}^{-1}L_{7} & -A_{7}^{-1} U_{7} & \dots\\
          0 &     0 &     0 &     I & \dots\\
          \dots & \dots & \dots & \dots & \dots\\
        \end{array}
        \right],
      \label{eqn:lineSolverProlongation}
    \end{aligned} \\
    \begin{aligned}
      R =& P^\intercal,
      \label{eqn:lineSolverRestriction}
    \end{aligned}
  \end{align}
  and $M_0^{-1} = \left(D - C A^{-1} B\right)^{-1}$ as in Proposition \ref{prop:linecoarsespace}
  converges in one iteration and thus is a direct solver.

  This method is a reformulation of the cyclic reduction algorithm
  of Buzbee, Golub and Nielson \cite{BuzbeeGolubNielson70}.
  \label{theo:twoleveldirect}
\end{theorem}
\begin{proof}
  The proof consists in applying Proposition \ref{prop:blockinverse} to
  $Q^\intercal M Q$ as in Proposition \ref{prop:linecoarsespace} and
  then using Lemma \ref{lem:inversionrichardsonequivalence}. We have
  \begin{align*}
    \left(Q^\intercal M Q\right)^{-1} =&
    \begin{aligned}
      \scalemath{0.7}{
        \left(
        \begin{array}{cc}
          A^{-1} & 0 \\
          0 & 0
        \end{array}
        \right)
        +
        \left(
        \begin{array}{c}
          -A^{-1} B \\
          I
        \end{array}
        \right)
        (D - C A^{-1} B)^{-1}
        \left(
        \begin{array}{cc}
          -A^{-1}B & I
        \end{array}
        \right)
        \left(
        I
        -
        Q^\intercal M Q
        \left(
        \begin{array}{cc}
          A^{-1} & 0 \\
          0 & 0
        \end{array}
        \right)
        \right)
      }
    \end{aligned} \\
    M^{-1} =&
    \begin{aligned}
      \scalemath{0.6}{
        Q \left(
        \begin{array}{cc}
          A^{-1} & 0 \\
          0 & 0
        \end{array}
        \right) Q^\intercal
        +
        Q \left(
        \begin{array}{c}
          -A^{-1} B \\
          I
        \end{array}
        \right)
        (D - C A^{-1} B)^{-1}
        \left(
        \begin{array}{cc}
          -A^{-1}B & I
        \end{array}
        \right)
        \left(
        I
        -
        Q^\intercal M Q
        \left(
        \begin{array}{cc}
          A^{-1} & 0 \\
          0 & 0
        \end{array}
        \right)
        \right) Q^\intercal
      }
    \end{aligned}\\
    M^{-1} =&
    \begin{aligned}
      \scalemath{0.6}{
        Q \left(
        \begin{array}{cc}
          A^{-1} & 0 \\
          0 & 0
        \end{array}
        \right) Q^\intercal
        +
        Q \left(
        \begin{array}{c}
          -A^{-1} B \\
          I
        \end{array}
        \right)
        (D - C A^{-1} B)^{-1}
        \left(
        \begin{array}{cc}
          -A^{-1}B & I
        \end{array}
        \right)
        Q^\intercal Q
        \left(
        I
        -
        Q^\intercal M Q
        \left(
        \begin{array}{cc}
          A^{-1} & 0 \\
          0 & 0
        \end{array}
        \right)
        \right) Q^\intercal
      }
    \end{aligned} \\
    M^{-1} =&
    \begin{aligned}
      \scalemath{0.6}{
        Q \left(
        \begin{array}{cc}
          A^{-1} & 0 \\
          0 & 0
        \end{array}
        \right) Q^\intercal
        +
        Q \left(
        \begin{array}{c}
          -A^{-1} B \\
          I
        \end{array}
        \right)
        (D - C A^{-1} B)^{-1}
        \left(
        \begin{array}{cc}
          -A^{-1}B & I
        \end{array}
        \right)
        Q^\intercal
        \left(
        I
        -
        M Q
        \left(
        \begin{array}{cc}
          A^{-1} & 0 \\
          0 & 0
        \end{array}
        \right) Q^\intercal
        \right) 
      }
    \end{aligned},
  \end{align*}
  from which it is easy to see that
  \begin{align*}
    S^{-1}=Q
    \left(
    \begin{array}{cc}
      A^{-1} & 0 \\
      0 & 0
    \end{array}
    \right) Q^\intercal
    \quad \text{ and } \quad
    P = Q \left(
        \begin{array}{c}
          -A^{-1} B \\
          I
        \end{array}
        \right),
  \end{align*}
  and by Lemma \ref{lem:inversionrichardsonequivalence}, a two-level
  preconditioned Richardson iteration without post-smoothing converges
  in one step and is therefore a direct solver.

  By Lemma \ref{lem:inversionrichardsonequivalence}, the two-level
  multigrid-preconditioned Richardson method without post-smoothing is a
  reformulation of \cite[Lemma 1]{GanderKwokZhang2018}, and by Lemma
  \ref{cor:inversionequivalence} it is equivalent to block-cyclic
  reduction.
  
\end{proof}

Corollary \ref{cor:multigridtridiagonal} concludes that the recursive
application of Theorem \ref{theo:twoleveldirect} gives a direct
multigrid solver.
\begin{cor}[Direct multigrid solver]
  The method described in Theorem \ref{theo:twoleveldirect} can be
  recursively applied to the coarse space, therefore providing a
  multigrid preconditioner that converges in one iteration.
  This is equivalent to the recursive application of two-level
  cyclic-reduction.
  \label{cor:multigridtridiagonal}
\end{cor}
\begin{proof}
  By Proposition
  \ref{prop:linecoarsespace}, the coarse space $M_0$ in Theorem
  \ref{theo:twoleveldirect} is block-tridiagonal and therefore
  algorithm \ref{alg:twolevelpreconditionernosmoothing} can be applied
  to solve it.
\end{proof}

We deduce the cost of the multigrid line solver to guarantee that it
is optimal in Proposition \ref{prop:lineSmootherCost}.
\begin{prop}
  Let $n$ be the total number of blocks of a block-tridiagonal system
  and $m$, the size of the blocks involved. Then the multigrid direct
  solver in Corollary \ref{cor:multigridtridiagonal} has
  $\mathcal{O}\left(nm^3\right)$ computational cost and
  $\mathcal{O}\left(n m^2\right)$ memory footprint when the operators
  at each level are built and stored.

  If every element is built on-the-fly, the algorithm has
  $\mathcal{O}\left(n^2 m^3\right)$ computational cost and
  $\mathcal{O}\left(n\right)$ memory footprint when implemented
  matrix-free.
  \label{prop:lineSmootherCost}
\end{prop}
\begin{proof}
  Let the blocks of a block-tridiagonal matrix be of size $m\times m$,
  and assume they are not sparse. Then the cost of applying each of
  them to a vector of size $m$ is $\mathcal{O}\left(m^2\right)$ and
  the cost of obtaining their inverse is $\mathcal{O}\left(m^3\right)$.

  The cost of building $S^{-1}$, $P$, $R$ and $M_0$ in Theorem
  \ref{theo:twoleveldirect} is dominated by the inverses $A_i^{-1}$
  and thus is $\mathcal{O}\left(\frac{n}2 m^3\right)$; this
  constitutes the cost of processing the finest mesh of the
  multigrid solver in Corollary \ref{cor:multigridtridiagonal}. The
  memory footprint of storing the finest mesh operators is
  $\mathcal{O}\left(\frac{n}2 m^2\right)$.

  When accounting for coarser levels of the multigrid method the costs
  of building the smoother, prolongation and restriction are
  \[\mathcal{O}\left(\frac{n}2 m^3\right) + \mathcal{O}\left(\frac{n}4
  m^3\right) + \mathcal{O}\left(\frac{n}8 m^3\right) + \dots =
  \mathcal{O}\left(n m^3\right).\] Subsequently, the cost of applying
  the operators is $\mathcal{O}\left(nm^3\right)$ and the memory
  footprint is $\mathcal{O}\left(nm^2\right)$.

  If we build every element on-the-fly, we have to build the sub-tree
  associated with every element in every mesh. There are
  $\mathcal{O}\left(n\right)$ elements considering all meshes; thus,
  the complexity of the algorithm is at most $\mathcal{O}\left(n^2
  m^3\right)$. The memory footprint is reduced to $\mathcal{O}(n)$
  given by the need to store the solution and the iterate.
\end{proof}

\begin{rem}
  When looking at the memory footprint of Proposition
  \ref{prop:lineSmootherCost}, the storage of $P$, $R$, $S^{-1}$ and
  $M_0$ for every level (in sparse format of course) multiplies $n$ by
  a constant at least equal to 6: 1 for each of $P$, $R$, $S^{-1}$ and
  3 for the tridiagonal $M_0$.

  This means that we can reduce the memory footprint at least by a
  factor $6$ by changing the complexity from
  $\mathcal{O}\left(nm^3\right)$ to
  $\mathcal{O}\left(n^2m^3\right)$. Such change might make sense for
  architectures with a high level of parallelism but low memory
  footprint, for certain values of $n$.
\end{rem}
This concludes the description of the line-solver that will be used to
invert the tridiagonal systems coming from the restriction of the 2D
finite element discretization of the Laplacian on the mesh given by
the GLL spectral discretization.

\subsection{Multigrid $\gamma$-cycle preconditioner}\label{sec:multigridPreconditioner}
We have all the elements needed to design our multigrid preconditioner
using the line preconditioners as line \emph{smoothers} (in the
context of multigrid), based on $A_p$ or $\tilde{A}_p$. We use a
$\gamma$-cycle instead of just a $V$-cycle since it improves
convergence properties without changing the order of complexity of the
overall method, as we will prove. In this section we will show how the
definitions below fall into place as pieces of a multilevel method. We
begin with the description of a hierarchy of meshes and the operators
needed to project and restrict functions between them.

Let $\mathbb{T}_\ell$ be a hierarchy of quadrilateral meshes in 2
dimensions. Considering multilevel methods, the index $\ell$ refers to
the mesh given by the tensor product of 1D GLL quadrature points of
polynomial degree $p_\ell$. Consider a sequence of subspaces
\[ V_1 \subset V_2 \subset \dots \subset V_{\ell} \subset \dots \subset V_{L},\]
given by different polynomial degree spaces as in
eq. \eqref{eqn:DiscreteSpace}. We introduce the projections
\begin{align*}
  R_\ell v_\ell: V_{\ell} \rightarrow V_{\ell-1} \quad \text{ and } \quad
  P_\ell v_\ell: V_{\ell-1} \rightarrow V_{\ell},
\end{align*}
defined simply as the interpolation operators from one basis to
another, at the corresponding quadrature nodes. The multi-dimensional
restriction operator is given in terms of the 1D operators in
tensor-product form and therefore has a complexity
$\mathcal{O}\left(p_\ell^3\right)$, and $R_\ell =
\left(P_\ell\right)^\intercal$ (see \cite[\S 3.3.1]{Ronquist1988} for
details).

Let $A_\ell$ be the operator defined in equation
\eqref{eqn:discreteBilinear} on the mesh $\mathbb{T}_\ell$. For the
rest of the paper we will redefine the operators
$P_{\mathscr{h}_i},P_{\mathscr{v}_i}$ and
$R_{\mathscr{h}_i},R_{\mathscr{h}_i}$ by adding an index $\ell$ such
that we can identify the corresponding multigrid level $\ell$ and line
$\mathscr{h}_i$ or $\mathscr{v}_i$.

\begin{figure}
  \begin{subfigure}[b]{\textwidth}
    \centering
    \resizebox{0.25\textwidth}{!}{
  \includegraphics{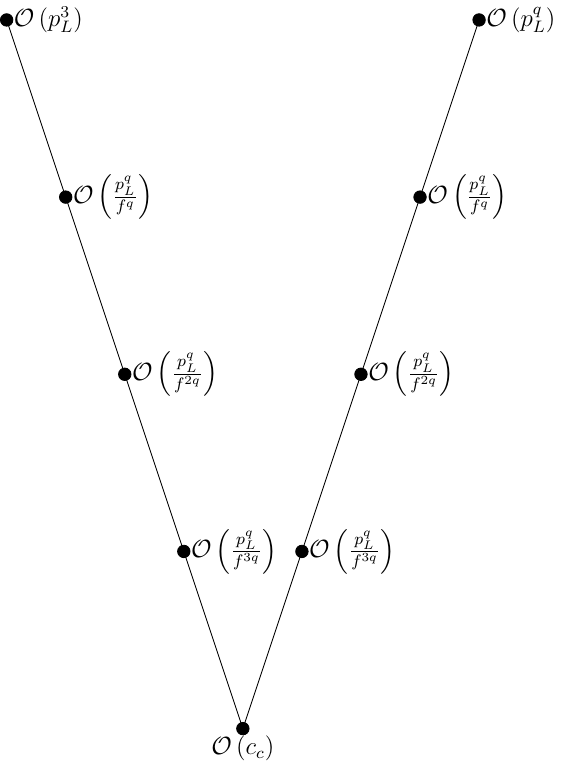}
    }
    \caption{$\gamma=1$}
  \end{subfigure}
  \begin{subfigure}[b]{\textwidth} 
    \centering
    \resizebox{\textwidth}{!}{
  \includegraphics{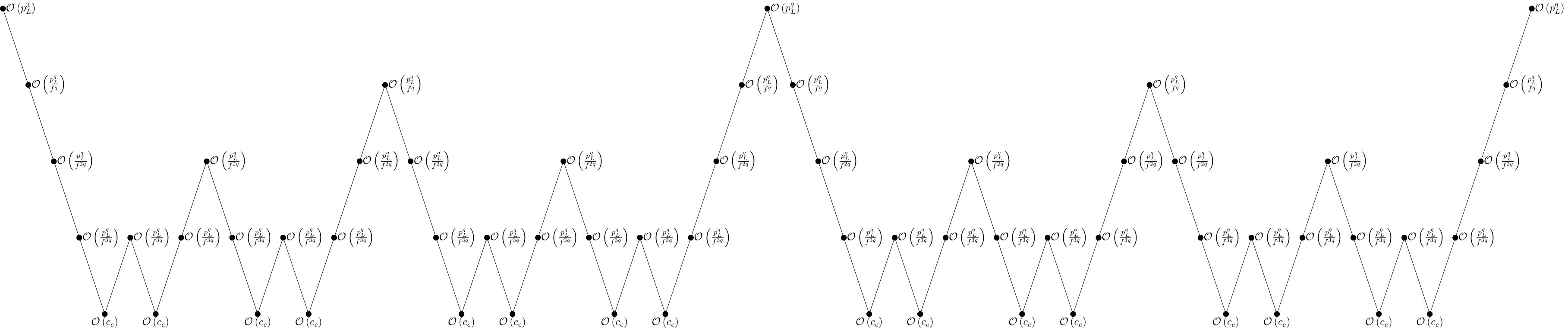}
    }
    \caption{$\gamma=2$}
  \end{subfigure}
  \begin{subfigure}[b]{\textwidth}
    \centering
    \resizebox{\textwidth}{!}{
  \includegraphics{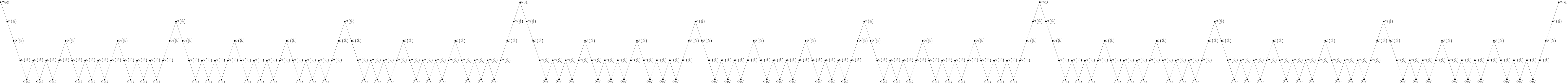}
    }
    \caption{$\gamma=3$}
  \end{subfigure}
  \caption{$5$-level multigrid $\gamma$-cycle for $\gamma=1,2,3$.}
  \label{fig:gammacycle}
\end{figure}

The multigrid $\gamma$-cycle is nowadays a standard technique in the
literature, therefore we only provide an inductive, compact definition
of the cycle below.
\begin{defi}[Multigrid $\gamma$-cycle (see Fig. \ref{fig:gammacycle})]
  Let $B_\ell$ be a smoother defined using
  the line preconditioners presented in \S\ref{sec:LinePreconditioners}, i.e.,
  \[ B_{\ell,\mathscr{h}}^{-1} = \sum_{i=2}^{p_\ell} P_{\ell,\mathscr{h}_i} \left(R_{\ell,\mathscr{h}_i} A_\ell P_{\ell,\mathscr{h}_i}\right)^{-1} R_{\ell,\mathscr{h}_i} \text{ and }
  B_{\ell,\mathscr{v}}^{-1} = \sum_{i=2}^{p_\ell} P_{\ell,\mathscr{v}_i} \left(R_{\ell,\mathscr{v}_i} A_\ell P_{\ell,\mathscr{v}_i}\right)^{-1} R_{\ell,\mathscr{v}_i}.
  \]
  If we use FEM line preconditioners let
  \[ \tilde{B}_{\ell,\mathscr{h}}^{-1} = \sum_{i=2}^{p_\ell} P_{\ell,\mathscr{h}_i} \left(R_{\ell,\mathscr{h}_i} \tilde{A}_\ell P_{\ell,\mathscr{h}_i}\right)^{-1} R_{\ell,\mathscr{h}_i} \text{ and }
  \tilde{B}_{\ell,\mathscr{v}}^{-1} = \sum_{i=2}^{p_\ell} P_{\ell,\mathscr{v}_i} \left(R_{\ell,\mathscr{v}_i} \tilde{A}_\ell P_{\ell,\mathscr{v}_i}\right)^{-1} R_{\ell,\mathscr{v}_i}
  \]
  Then Alg. \ref{alg:gammacycle} is the multigrid $\gamma$-cycle.
  \textnormal{
  \begin{algorithm}[H]
    \caption{Multigrid $\gamma$-cycle (see
      Fig. \ref{fig:gammacycle}). \\
      Let the multigrid preconditioner $M_L^{-1}$ be defined by
      induction, let $\alpha \in \mathbb{R}^+$ be a relaxation
      parameter and $M_0^{-1}=A_0^{-1}$. For $0\le \ell \le L$, we
      define the action $M_\ell^{-1} r$ of $M_\ell^{-1}$ on a vector $r\in
      V_\ell$ in terms of $M_{\ell-1}^{-1}$ below. \\
      \textbf{Input:} $r\in V_\ell$. \textbf{Output:} $y\in V_\ell$.}
    \begin{algorithmic}[1] % [1] shows line numbers
      \State Let $x_0 = 0$.
      \State Create or assign $x_i$ for $i=1,\dots,2m$, by 
      \Statex $m$ horizontal pre-smoothing steps using $B_{\ell,\mathscr{h}}^{-1}$
      or $\tilde{B}_{\ell,\mathscr{h}}^{-1}$
      \[x_i=x_{i-1}+\alpha B_{\ell,\mathscr{h}}^{-1}\left(r - A_\ell x_{i-1}\right),\]
      \Statex and $m$ vertical pre-smoothing steps using $B_{\ell,\mathscr{v}}^{-1}$
      or $\tilde{B}_{\ell,\mathscr{v}}^{-1}$
      \[x_i=x_{i-1}+\alpha B_{\ell,\mathscr{v}}^{-1}\left(r - A_\ell x_{i-1}\right).\]
      \For{$g=1$ to $\gamma$}
      \State Create or assign $y_0$ by correcting the residual with a coarser grid
      \Statex \quad\quad (i.e. here the algorithm calls itself unless
      $\ell=1$ since $M_0^{-1}=A_0^{-1}$)
      \[y_0 = x_{2m} + P_\ell M_{\ell-1}^{-1} R_\ell\left(r - A_\ell x_{2m}\right).\]
      \State Create or assign $y_i$ for $i=1,\dots,2m$, by 
      \Statex \hspace{0.5cm}$m$ vertical post-smoothing steps using $B_{\ell,\mathscr{v}}^{-1}$
      or $\tilde{B}_{\ell,\mathscr{v}}^{-1}$
      \[y_i=y_{i-1}+\alpha B_{\ell,\mathscr{v}}^{-1}\left(r - A_\ell y_{i-1}\right),\]
      \Statex \hspace{0.5cm}and $m$ horizontal post-smoothing steps using $B_{\ell,\mathscr{h}}^{-1}$
      or $\tilde{B}_{\ell,\mathscr{h}}^{-1}$
      \[y_i=y_{i-1}+\alpha B_{\ell,\mathscr{h}}^{-1}\left(r - A_\ell y_{i-1}\right).\]
      \EndFor
      \State $y = y_{2m}$.
    \end{algorithmic}
    \label{alg:gammacycle}
  \end{algorithm}}
  \label{def:gammacycle}
\end{defi}

Figure \ref{fig:gammacycle} shows a graphical depiction of the
multigrid $\gamma$-cycle as in Def. \ref{def:gammacycle} for a 5-level
structure.  The next section will dive into the computational
complexity of the method, ultimately proving that it is of the same
order as the application of $A_p$.

\subsubsection{Computational complexity}
We deduce in this section a proposition showing how the computational
complexity will depend on $\gamma$, in order to choose its optimal
value.
\begin{prop}
  Let a $\gamma$-multigrid cycle as defined in
  \S\ref{sec:multigridPreconditioner} be such that the cost of
  restriction, prolongation and smoothing at an arbitrary level $\ell$
  has a complexity $\mathcal{O}\left(p_\ell^q\right)$ for $q \in
  \mathbb{N}$, $f = \frac{p_\ell}{p_{\ell-1}}$ be the
  coarsening factor, and let the cost of solving the coarsest level be
  independent of $L$.
  The complexity of the $\gamma$-cycle, for a
  constant cost of solving the coarsest level $\ell=0$, is
  \begin{equation}
    \begin{cases}
      \text{if } 1 < \gamma < f^q \quad & \mathcal{O}\left(\frac{\gamma^2}{f^q} p_L^q\right) \\
      \text{if } \gamma = f^q \quad & \mathcal{O}\left(\frac{\gamma^2}{f^q}p_L^q\log_{f}{p_L} \right) \\
      \text{if } \gamma > f^q \quad &
      \mathcal{O}\left(\gamma p_L^q p_L^{\log_f\left(\frac{\gamma}{f^q}\right)}\right) \\
    \end{cases}
  \end{equation}
  \label{prop:multigridComplexity}
\end{prop}
\begin{proof}
  Figure \ref{fig:gammacycle} shows that the cost of an arbitrary
  level $\ell$, except the coarsest since it is assumed constant, is
  \[
  \mathcal{O}\left(\gamma^{\ell} (\gamma+1) \frac{p_L^q}{f^{q\ell}}\right) =
  \mathcal{O}\left((\gamma+1)p_L^q \left(\frac{\gamma}{f^q}\right)^\ell\right).
  \]
  Adding over all levels and using the geometric sum formula we obtain
  \[
  \mathcal{O}\left((\gamma+1)p_L^q\sum_{\ell=0}^L \left(\frac{\gamma}{f^q}\right)^\ell\right) =
  \begin{cases}
    \text{if } 1 < \gamma < f^q \quad & \mathcal{O}\left((\gamma+1) \frac{\gamma}{f^q} p_L^q \right), \\
    \text{if } \gamma = f^q \quad & \mathcal{O}\left((\gamma+1) \frac{\gamma}{f^q} p_L^q L \right), \\
    \text{else } & \mathcal{O}\left((\gamma+1)p_L^q \frac{\left(\frac{\gamma}{f^q}\right)^L-1}{\frac{\gamma}{f^q}-1} \right).
  \end{cases}
  \]
  By definition of $f$, for constant $p_0$ we have $p_L = \mathcal{O}\left(f^L\right)$,
  thus $L=\mathcal{O}\left(\log_{f}{p_L}\right)$ and we get
  \[
  \mathcal{O}\left((\gamma+1)p_L^q\sum_{\ell=0}^L \left(\frac{\gamma}{f^q}\right)^\ell\right) =
  \begin{cases}
    \text{if } 1 < \gamma < f^q \quad & \mathcal{O}\left(\frac{\gamma^2}{f^q} p_L^q\right), \\
    \text{if } \gamma = f^q \quad & \mathcal{O}\left(\frac{\gamma^2}{f^q}p_L^q\log_{f}{p_L} \right), \\
    \text{if } \gamma > f^q \quad &
    \mathcal{O}\left(\gamma p_L^q p_L^{\log_f\left(\frac{\gamma}{f^q}\right)}\right), \\
  \end{cases}
  \]
  since $\left(\frac{\gamma}{f^q}\right)^{\log_f\left(p_L\right)} =
  p_L^{\log_f\left(\frac{\gamma}{f^q}\right)}$, and the result is
  achieved.
\end{proof}

Proposition \ref{prop:multigridComplexity} shows that in 2D, for the
standard choice of $f=2$ we can choose $\gamma$ as large as $7$ and our
solver will remain of complexity $\mathcal{O}\left(p_L^3\right)$ like
the application of $A_L$! Should we choose $\gamma=8$ we would only be
adding a logarithmic order to the solver.

All the structures required for the application of the solver only
require the storage of 1D structures, thus keeping memory the footprint
controlled.

\section{Numerical experiments}\label{sec:NumericalExperiments}

\begin{table}[h!]
  \centering
  \begin{tabular}{>{\centering\arraybackslash}m{0.17\textwidth} >{\centering\arraybackslash}m{0.37\textwidth} >{\centering\arraybackslash}m{0.37\textwidth} }  
    Problem & GLL smoother & FEM smoother \\
    \hline
    \centering
    \begin{minipage}[c][2.5cm][c]{\linewidth}  
      \centering
      \includegraphics[width=0.7\linewidth,keepaspectratio]{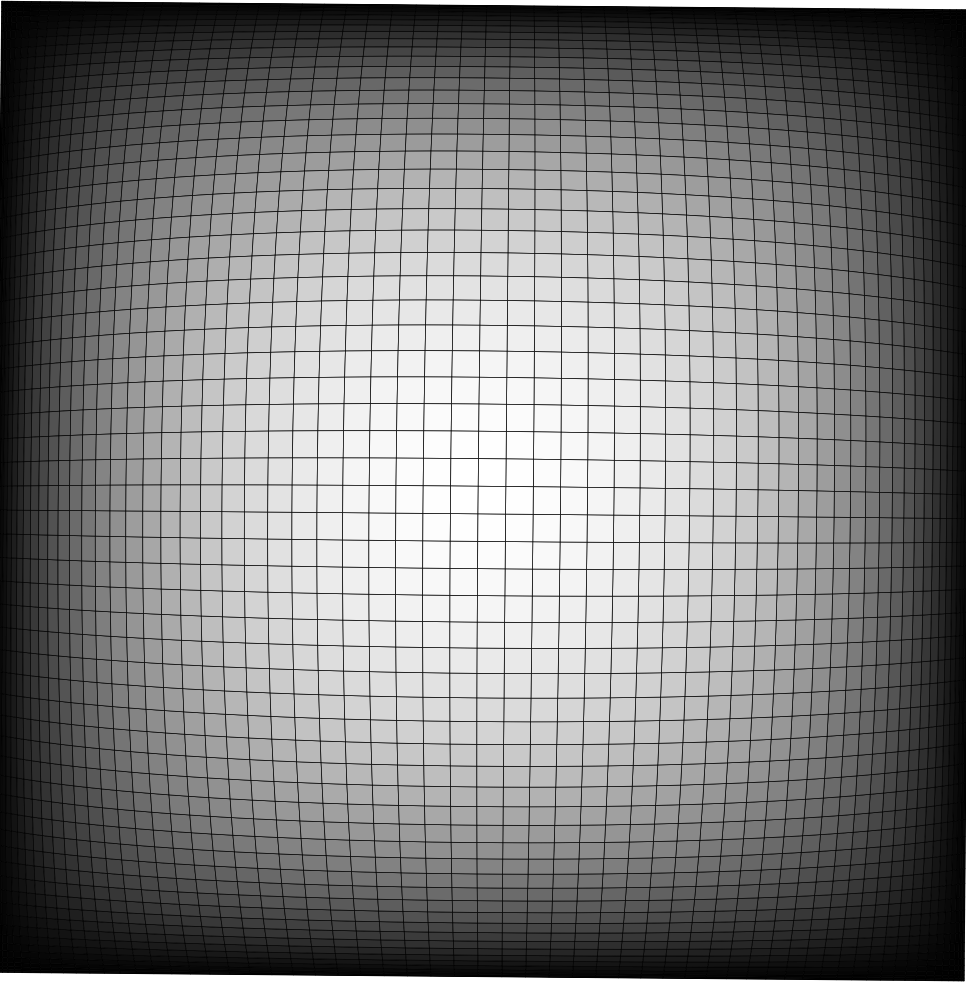}
      \[
      \scalemath{0.8}{\Delta u = -1}
      \]
    \end{minipage}
    &
    \begin{minipage}{\linewidth}
      \resizebox{0.99\linewidth}{!}{  
        \begin{tabular}{c|c c c c c c c c}
          \diagbox{$p_L$}{$\gamma$} & 1 & 2 & 3 & 4 & 5 & 6 & 7 & 8 \\ \hline
          8  &  6 &  5 &  4 & 4 & 3 & 3 & 3 & 3 \\
          16 & 11 &  8 &  7 & 6 & 5 & 5 & 4 & 4 \\
          32 & 19 & 12 &  9 & 7 & 6 & 5 & 5 & 5 \\
          64 & 31 & 17 & 11 & 8 & 7 & 6 & 5 & 5 \\
        \end{tabular}
      }
    \end{minipage}
    &
    \begin{minipage}{\linewidth}
      \resizebox{\linewidth}{!}{  
        \begin{tabular}{c|c c c c c c c c}
          \diagbox{$p_L$}{$\gamma$} & 1 & 2 & 3 & 4 & 5 & 6 & 7 & 8 \\ \hline
          8  &  9 &  7 &  6 &  5 & 5 & 5 & 4 & 4 \\
          16 & 14 & 10 &  8 &  7 & 6 & 5 & 5 & 4 \\
          32 & 23 & 14 & 10 &  8 & 7 & 6 & 5 & 5 \\
          64 & 40 & 20 & 13 &  9 & 7 & 6 & 5 & 5 \\
        \end{tabular}
      }
    \end{minipage}
    \\
    \begin{minipage}[c][2.5cm][c]{\linewidth}  
      \centering
      \includegraphics[width=0.7\linewidth,keepaspectratio]{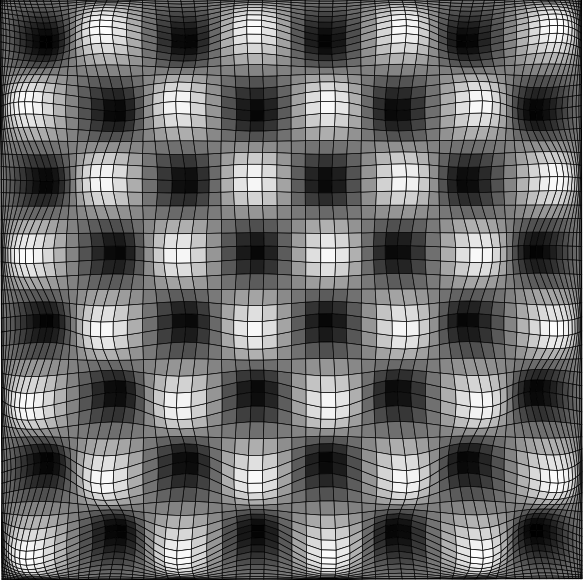}
      \begin{align*}
        \scalemath{0.6}{
          \begin{aligned}
            &\Delta u = \\
            &\Delta\left(\sin\left(8k\pi x\right)\sin\left(8k\pi y\right)\right)
          \end{aligned}
          }
      \end{align*}
    \end{minipage}
    &
    \begin{minipage}{\linewidth}
      \resizebox{0.99\linewidth}{!}{  
        \begin{tabular}{c|c c c c c c c c}
          \diagbox{$p_L$}{$\gamma$} & 1 & 2 & 3 & 4 & 5 & 6 & 7 & 8 \\ \hline
          8  &  6 &  5 &  4 & 4 & 3 & 3 & 3 & 3 \\
          16 & 11 &  8 &  7 & 6 & 5 & 5 & 5 & 4 \\
          32 & 17 & 12 &  9 & 8 & 7 & 6 & 6 & 5 \\
          64 & 27 & 16 & 11 & 9 & 8 & 7 & 6 & 5 \\
        \end{tabular}
      }
    \end{minipage}
    &
    \begin{minipage}{\linewidth}
      \resizebox{\linewidth}{!}{  
        \begin{tabular}{c|c c c c c c c c}
          \diagbox{$p_L$}{$\gamma$} & 1 & 2 & 3 & 4 & 5 & 6 & 7 & 8 \\ \hline
          8  &  8 &  6 &  6 &  5 & 5 & 4 & 4 & 4 \\
          16 & 13 & 10 &  8 &  7 & 6 & 6 & 5 & 5 \\
          32 & 20 & 13 & 10 &  8 & 7 & 6 & 6 & 5 \\
          64 & 33 & 19 & 13 & 10 & 8 & 7 & 6 & 5 \\
        \end{tabular}
      }
    \end{minipage}
    \\
    \begin{minipage}[c][2.5cm][c]{\linewidth}  
      \centering
      \includegraphics[width=0.7\linewidth,keepaspectratio]{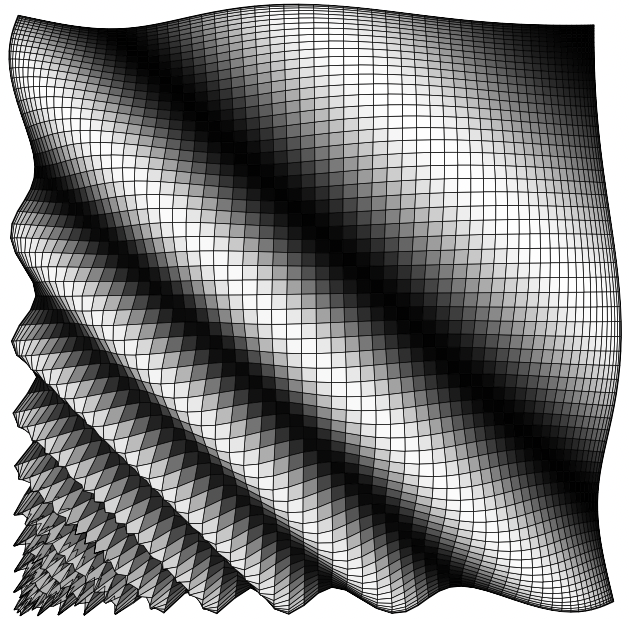}
      \begin{align*}
        \scalemath{0.6}{
          \begin{aligned}
            &\Delta u = \\
            &\Delta\left(\sin\left(\frac{8\pi}{x+y+\frac{\pi}{10}}\right)\right)
          \end{aligned}
          }
      \end{align*}
    \end{minipage}
    &
    \begin{minipage}{\linewidth}
      \resizebox{0.99\linewidth}{!}{  
        \begin{tabular}{c|c c c c c c c c}
          \diagbox{$p_L$}{$\gamma$} & 1 & 2 & 3 & 4 & 5 & 6 & 7 & 8 \\ \hline
          8  & 10 &  7 &  6 &  5 &  5 & 4 & 4 & 4 \\
          11 & 16 & 11 &  9 &  7 &  6 & 6 & 5 & 5 \\
          32 & 27 & 17 & 12 & 10 &  8 & 7 & 6 & 6 \\
          64 & 45 & 24 & 15 & 12 & 10 & 9 & 9 & 8 \\
        \end{tabular}
      }
    \end{minipage}
    &
    \begin{minipage}{\linewidth}
      \resizebox{\linewidth}{!}{  
        \begin{tabular}{c|c c c c c c c c}
          \diagbox{$p_L$}{$\gamma$} & 1 & 2 & 3 & 4 & 5 & 6 & 7 & 8 \\ \hline
          8  & 13 & 10 &  8 &  7 &  6 & 6 & 5 & 5 \\
          11 & 20 & 13 & 10 &  8 &  7 & 6 & 6 & 5 \\
          32 & 32 & 19 & 13 & 11 &  9 & 8 & 7 & 6 \\
          64 & 56 & 28 & 18 & 13 & 10 & 8 & 7 & 6 \\
        \end{tabular}
      }
    \end{minipage}
  \end{tabular}
  \caption{{\bf Square domain.} Iteration count to reduce the residual
    by a factor $10^{-8}$ for a domain $\Omega = [0,1]^2$ without
    deformation. We use $m = 1$ and $\alpha_{GLL}=\frac23$
    and $\alpha_{FEM}=0.16$ for all experiments.}
  \label{tab:numericalExperiments1}
\end{table}
Table \ref{tab:numericalExperiments1} illustrates the behavior of the
method for different right hand sides on the domain $[0,1]^2$. We show
right hand sides generated by constant values and manufactured
solutions, as well as homogeneous and non-homogeneous boundary
conditions. For all cases, our convergence criteria is the reduction
of the residual by a factor $10^{-8}$. We observe a significant reduction of
the iteration count, in particular for $\gamma = 7$ that, as shown in
Proposition \ref{prop:multigridComplexity}, keeps the complexity of the
preconditioned operator equal to that of the unpreconditioned one.

We do not intend to optimize the relaxation parameters $\alpha$, so we
choose them by simple observation of the convergence behavior to be
$\alpha_{GLL}=\frac23$ and $\alpha_{FEM}=0.16$ for all experiments. A
proper optimization of these values would require a dedicated study as
in \cite{LuceroLorcaGander2024,GanderLucero22_2}.

We observe a noticeable increase in the overall iteration count,
particularly for low gamma values in the last case. This occurs
because we begin the solver with an initial guess of zero and thus the
first two cases begin without high frequencies in the residual. The
lack of high frequencies requires a smaller Krylov subspace to
eliminate them. Experiments with random initial guesses show that the
iteration counts in all cases are similar to the third case. Our
examples show that iteration counts diminish when the solution is
smooth (for a zero guess), which will often be the case when we
implement our solver for multiple spectral elements.

\begin{table}[h!]
  \centering
  \begin{tabular}{>{\centering\arraybackslash}m{0.45\textwidth} | >{\centering\arraybackslash}m{0.08\textwidth} |
      >{\centering\arraybackslash}m{0.11\textwidth} |
      >{\centering\arraybackslash}m{0.11\textwidth} | >{\centering\arraybackslash}m{0.1\textwidth}}
    \multirow{2}{*}{Problem} & Angle & GLL smoother & FEM smoother & Relative\\
    &            $[^{\circ}C]$   & iterations   & iterations   & error \\ \hline
    \multirow{16}{*}{
    \centering
    \begin{minipage}[c][2.5cm][c]{\linewidth}  
      \centering
      \includegraphics[width=0.4\linewidth,keepaspectratio]{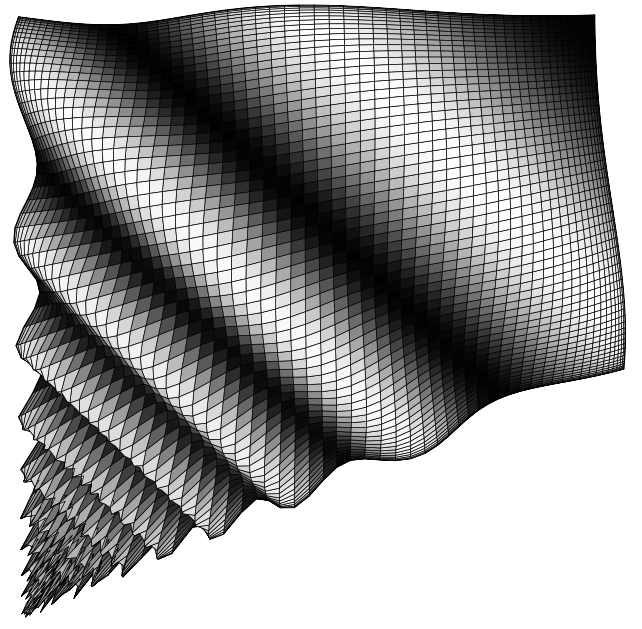}
      \begin{align*}
        \scalemath{1}{
          \begin{aligned}
            &\Delta u = &\Delta\left(\sin\left(\frac{8\pi}{x+y+\frac{\pi}{10}}\right)\right)
          \end{aligned}
        }
      \end{align*}
    \end{minipage}}
    &       0 &   9 &    7 &$2\cdot 10^{-13}$ \\
    &      10 &   9 &    7 &$2\cdot 10^{-10}$ \\
    &      11 &   9 &    7 &$3\cdot 10^{-10}$ \\
    &      12 &  10 &    7 &$5\cdot 10^{-10}$ \\
    &      13 &  11 &    7 &$9\cdot 10^{-10}$ \\
    &      14 &  15 &    7 &$2\cdot 10^{-09}$ \\
    &      15 &  24 &    7 &$3\cdot 10^{-09}$ \\
    &      16 & >30 &    7 &$6\cdot 10^{-09}$ \\
    &      17 & >30 &    7 &$1\cdot 10^{-08}$ \\
    &      18 & >30 &    9 &$2\cdot 10^{-08}$ \\
    &      19 & >30 &   11 &$3\cdot 10^{-08}$ \\
    &      20 & >30 &   14 &$5\cdot 10^{-08}$ \\
    &      21 & >30 &   17 &$8\cdot 10^{-08}$ \\
    &      22 & >30 &   20 &$1\cdot 10^{-07}$ \\
    &      23 & >30 &  >30 &$1\cdot 10^{-07}$ \\\cline{1-2}
    \multirow{2}{*}{Problem} & \multirow{2}{*}{Height} & \multirow{2}{*}{} & \multirow{2}{*}{} & \multirow{2}{*}{} \\
    &              & \multirow{2}{*}{} & \multirow{2}{*}{} & \multirow{2}{*}{}  \\ \cline{1-2}
    \multirow{8}{*}{
    \centering
    \begin{minipage}[c][2.5cm][c]{\linewidth}  
      \centering
      \includegraphics[width=0.4\linewidth,keepaspectratio]{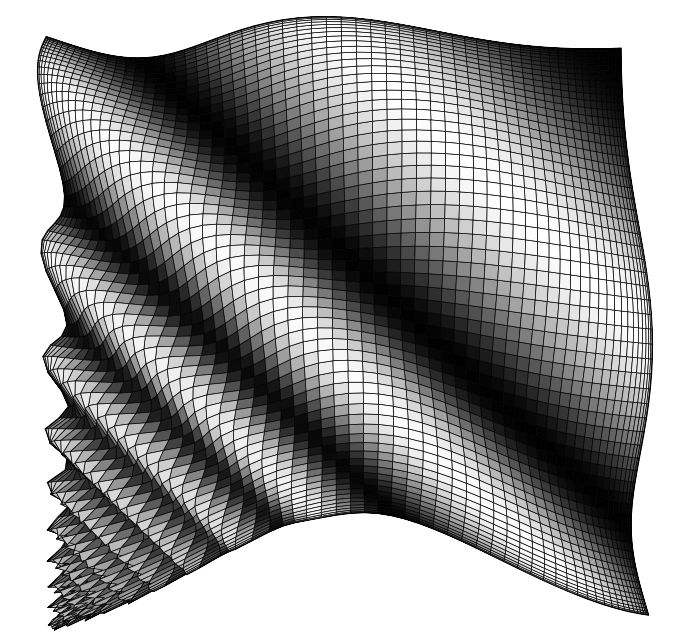}
      \begin{align*}
        \scalemath{1.}{
          \begin{aligned}
            &\Delta u = &\Delta\left(\sin\left(\frac{8\pi}{x+y+\frac{\pi}{10}}\right)\right)
          \end{aligned}
        }
      \end{align*}
    \end{minipage}}
    &       0 &   9 &   7 &$2 \cdot 10^{-13}$ \\
    &    0.10 &   9 &   7 &$2 \cdot 10^{-10}$ \\
    &    0.15 &   9 &   9 &$4 \cdot 10^{-9}$ \\
    &    0.16 &   9 &   9 &$6 \cdot 10^{-9}$ \\
    &    0.17 &  11 &  10 &$8 \cdot 10^{-9}$ \\
    &    0.18 &  17 &  10 &$1 \cdot 10^{-8}$ \\
    &    0.19 &  21 & >30 &$1 \cdot 10^{-8}$ \\
    &    0.20 & >30 & >30 &$2 \cdot 10^{-8}$ \\
  \end{tabular}
  \caption{{\bf Deformed domain.} Iteration count to reduce the
    residual by a factor $10^{-8}$ for a domain $\Omega = [0,1]^2$ that
    has been deformed. We use $p_L=64$, $\gamma = 7$, $m = 1$ and
    $\alpha_{GLL}=\frac23$ and $\alpha_{FEM}=0.16$ for all
    experiments.}
  \label{tab:numericalExperiments2}
\end{table}
Table \ref{tab:numericalExperiments2} shows results for deformed
domains using $p_L=64$ and $\gamma=7$, where we increase the
deformation until the restriction, prolongation and Jacobians struggle
to represent the problem.  We observe that for moderate deformations
the scaling of the method is kept, and this motivates the use of SEM
for severely deformed domains, to keep the deformation bounded
\emph{on each element}.

We focus our attention on the behavior for $\gamma=7$, that provides
the desired complexity as we proved in
Prop. \ref{prop:multigridComplexity}, even though as we saw, using
$\gamma=8$ would only add a logarithmic term to the complexity.

The maximum polynomial degree is dictated by the capacity of our
algorithms to evaluate residuals for very high-degree polynomials, we
consider only a moderately high-degree polynomial of $p_L=64$ to be
able to continue to use double precision floating point
numbers. Higher degree polynomials would require higher precision and
special care on the truncation error induced when evaluating
residuals. We want to ultimately use our preconditioner as a smoother
for SEMs, so $p_L=64$ is more than enough for our applications.

Results look satisfactory and the iteration count seem to approach an
upper bound at a handful of $\gamma$-cycles. We illustrate the limits
of our method in Table \ref{tab:numericalExperiments2}. The
deformation has to be kept moderate, otherwise the iterations seem not
to have an upper bound but it is expected that the preconditioner will
not be able to cope with deformations that approach unbounded
Jacobians. This is expected since our restriction and prolongations to
coarser levels only approximate the geometry, and constitute a few
\emph{variational crimes} \cite{Strang1972, Strang1973,
  MadayRonquist1990} and \cite[Ch.10]{BrennerScott2008}.

Overall, our experiments show that the parameters of the multigrid
method need to be adapted to the total complexity, since even though
the \emph{scaling} tells us that higher values of $\gamma$ are
$\mathcal{O}\left(p^3\right)$, the value of $p$ has to be high enough
for the scaling to pay off.

\section{Conclusion}
We studied the application of a $p$-multigrid $\gamma$-cycle
preconditioner with line smoothers to solve a 2D
Gauss-Legendre-Lobatto spectral discretization of the Poisson
equation. We showed that the iteration count seems to have an upper
bound for moderate deformations. We illustrated our findings with a
variety of numerical experiments including non-homogeneous boundary
conditions and deformed geometries. The solver shows a satisfactory
performance at relatively high polynomial degrees in a single 2D
spectral element, with $\mathcal{O}\left(p^3\right)$ complexity. Which
shows its potential as a smoother for our ongoing development of an
$hp$-multigrid solver for a Spectral Element Method discretization of
the Poisson problem. Our next steps involve evaluating multi-element
discretizations in two and three dimensions.

\section*{Acknowledgments}
\thispagestyle{empty}

JPLL and DR wish to acknowledge gratefully our support by NOAA award number NA24OARX405C0008.

\clearpage

\bibliographystyle{abbrv}
\bibliography{pmultigrid}

\begin{thebibliography}{10}

\bibitem{Axelsson1985}
O.~Axelsson.
\newblock Incomplete block matrix factorization preconditioning methods. {T}he
  ultimate answer?
\newblock {\em Journal of Computational and Applied Mathematics}, 12-13:3--25,
  1985.

\bibitem{Axelsson1986}
O.~Axelsson.
\newblock Analysis of incomplete matrix factorizations as multigrid smoothers
  for vector and parallel computers.
\newblock {\em Applied Mathematics and Computation}, 19(1-4):3--22, 1986.

\bibitem{Bernstein2009}
D.~Bernstein.
\newblock {\em Matrix Mathematics: Theory, Facts, and Formulas - Second
  Edition}.
\newblock Princeton University Press, 2009.

\bibitem{Brandt1977}
A.~Brandt.
\newblock Multi-level adaptive solutions to boundary-value problems.
\newblock {\em Mathematics of computation}, 31(138):333--390, 1977.

\bibitem{BrennerScott2008}
S.~C. Brenner and L.~R. Scott.
\newblock {\em The Mathematical Theory of Finite Element Methods}, volume~15 of
  {\em Texts in Applied Mathematics}.
\newblock Springer, New York, 3rd edition, 2008.

\bibitem{BuzbeeGolubNielson70}
B.~L. Buzbee, G.~H. Golub, and C.~W. Nielson.
\newblock On direct methods for solving {P}oisson's equations.
\newblock {\em SIAM Journal on Numerical Analysis}, 7(4):627--656, 1970.

\bibitem{CanutoYousuffQuarteroniZang1987}
C.~Canuto, M.~Y. Hussaini, A.~Quarteroni, and T.~A. Zang.
\newblock {\em Spectral Methods in Fluid Dynamics}.
\newblock Springer-Verlag, Berlin, Heidelberg, 1987.

\bibitem{Demko1984}
S.~Demko, W.~F. Moss, and P.~W. Smith.
\newblock Decay rates for inverses of band matrices.
\newblock {\em Mathematics of Computation}, 43(168):491--499, 1984.

\bibitem{DevilleMund1985}
M.~Deville and E.~Mund.
\newblock Chebyshev pseudospectral solution of second-order elliptic equations
  with finite element preconditioning.
\newblock {\em Journal of Computational Physics}, 60(3):517--533, 1985.

\bibitem{DevilleFischerMund2002}
M.~O. Deville, P.~F. Fischer, and E.~H. Mund.
\newblock {\em High-Order Methods for Incompressible Fluid Flow}.
\newblock Cambridge Monographs on Applied and Computational Mathematics.
  Cambridge University Press, 2002.

\bibitem{Fischer1997}
P.~F. Fischer.
\newblock An overlapping {S}chwarz method for spectral element solution of the
  incompressible {N}avier--{S}tokes equations.
\newblock {\em Journal of Computational Physics}, 133(1):84--101, 1997.

\bibitem{GanderKwokZhang2018}
M.~J. Gander, F.~Kwok, and H.~Zhang.
\newblock Multigrid interpretations of the parareal algorithm leading to an
  overlapping variant and {MGRIT}.
\newblock {\em Comput. Visual Sci.}, 19:59--74, Jul 2018.

\bibitem{GottliebOrszag1977}
D.~Gottlieb and S.~A. Orszag.
\newblock {\em Numerical Analysis of Spectral Methods}.
\newblock Society for Industrial and Applied Mathematics, 1977.

\bibitem{Hackbusch1985}
W.~Hackbusch.
\newblock {\em Multi-Grid Methods and Applications}, volume~4 of {\em Springer
  Series in Computational Mathematics}.
\newblock Springer-Verlag, Berlin, Heidelberg, 1985.

\bibitem{Heinrichs1988}
W.~Heinrichs.
\newblock Line relaxation for spectral multigrid methods.
\newblock {\em Journal of Computational Physics}, 77(1):166--182, 1988.

\bibitem{Hemker2003}
P.~Hemker, W.~Hoffmann, and M.~van Raalte.
\newblock Two-level {F}ourier analysis of a multigrid approach for
  discontinuous {G}alerkin discretization.
\newblock {\em SIAM Journal on Scientific Computing}, 3(25):1018--1041, 2003.

\bibitem{Hemker2004}
P.~W. Hemker, W.~Hoffmann, and M.~H. van Raalte.
\newblock Fourier two-level analysis for discontinuous {G}alerkin
  discretization with linear elements.
\newblock {\em Numerical Linear Algebra with Applications}, 5 --
  6(11):473--491, 2004.

\bibitem{Lagrange1795}
J.-L. Lagrange.
\newblock {\em Leçons Elémentaires sur les Mathématiques: Leçon Cinquième.
  Sur l'usage des courbes dans la solution des problèmes. Republished in
  Serret, Joseph-Alfred, ed. (1877). Oeuvres de Lagrange.}, volume Vol 7.
\newblock Gauthier-Villars, 1795.

\bibitem{Lanczos1956}
C.~Lanczos.
\newblock {\em Applied Analysis}.
\newblock Applied Analysis. Prentice-Hall, 1956.

\bibitem{GanderLucero22_2}
J.~P. {Lucero Lorca} and M.~J. Gander.
\newblock Should multilevel methods for discontinuous {G}alerkin
  discretizations use discontinuous interpolation operators?
\newblock {\em Lecture Notes in Computational Science and Engineering}, 145,
  Mar 2022.

\bibitem{LuceroLorcaGander2024}
{{Lucero Lorca}, José Pablo} and {Gander, Martin Jakob}.
\newblock Optimization of two-level methods for {DG} discretizations of
  reaction-diffusion equations.
\newblock {\em ESAIM: M2AN}, 58(6):2351--2386, 2024.

\bibitem{MadayRonquist1990}
Y.~Maday and E.~M. Rønquist.
\newblock Optimal error analysis of spectral methods with emphasis on
  non-constant coefficients and deformed geometries.
\newblock {\em Computational Methods in Applied Mechanics and Engineering},
  80:91--115, 1990.

\bibitem{Olson2007}
L.~Olson.
\newblock Algebraic multigrid preconditioning of high-order spectral elements
  for elliptic problems on a simplicial mesh.
\newblock {\em SIAM Journal on Scientific Computing}, 29(5):2189--2209, 2007.

\bibitem{Orszag1980}
S.~A. Orszag.
\newblock Spectral methods for problems in complex geometries.
\newblock {\em Journal of Computational Physics}, 37(1):70--92, 1980.

\bibitem{PartnerRothman1995}
S.~V. Parter and E.~E. Rothman.
\newblock Preconditioning {L}egendre spectral collocation approximations to
  elliptic problems.
\newblock {\em SIAM Journal on Numerical Analysis}, 32(2):333--385, 1995.

\bibitem{Passot1987}
T.~Passot and A.~Pouquet.
\newblock Numerical simulation of compressible homogeneous flows in the
  turbulent regime.
\newblock {\em J. Fluid Mech.}, 181:441--466, 1987.

\bibitem{Patera1984}
A.~Patera.
\newblock A spectral element method for fluid dynamics: Laminar flow in a
  channel expansion.
\newblock {\em Journal of Computational Physics}, 54:468--488, 1984.

\bibitem{PavarinoWidlund1996}
L.~F. Pavarino and O.~B. Widlund.
\newblock A polylogarithmic bound for an iterative substructuring method for
  spectral elements in three dimensions.
\newblock {\em SIAM journal on numerical analysis}, 33(4):1303--1335, 1996.

\bibitem{Phillips1987}
T.~N. Phillips.
\newblock Relaxation schemes for spectral multigrid methods.
\newblock {\em Journal of Computational and Applied Mathematics},
  18(2):149--162, 1987.

\bibitem{Pope2000}
S.~Pope.
\newblock {\em Turbulent Flows}.
\newblock Cambridge University Press, 2000.

\bibitem{RonquistPatera1987_2}
E.~R{\o}nquist and A.~Patera.
\newblock A {L}egendre spectral element method for the {S}tefan problem.
\newblock {\em International Journal for Numerical Methods in Engineering},
  24:2273--2299, 1987.

\bibitem{Ronquist1988}
E.~M. R{\o}nquist.
\newblock {\em Optimal spectral element methods for the unsteady
  three-dimensional incompressible Navier-Stokes equations}.
\newblock PhD thesis, Massachusetts Institute of Technology, Cambridge, MA,
  1988.
\newblock Thesis (Ph.D.)--Massachusetts Institute of Technology, Dept. of
  Mechanical Engineering.

\bibitem{RonquistPatera1987}
E.~M. R{\o}nquist and A.~T. Patera.
\newblock Spectral element multigrid. {I}. formulation and numerical results.
\newblock {\em Journal of Scientific Computing}, 2(4):389--406, 1987.

\bibitem{Rosenberg2023}
D.~Rosenberg, B.~Flynt, M.~Govett, and I.~Jankov.
\newblock Geofluid object workbench ({G}eo{F}low) for atmospheric dynamics in
  the approach to exascale: Spectral element formulation and {CPU} performance.
\newblock {\em Month. Weather Rev.}, 151:2521--2540, 2023.

\bibitem{Roth1934}
W.~E. Roth.
\newblock {On direct product matrices}.
\newblock {\em Bulletin of the American Mathematical Society}, 40(6):461 --
  468, 1934.

\bibitem{Legendre1785}
A.~royale des~sciences (France).
\newblock {\em Mémoires de mathématique et de physique, presentés à
  l'Académie royale des sciences, par divers sçavans \& lûs dans ses
  assemblées}, volume Vol 10.
\newblock Paris,, 1785.
\newblock https://www.biodiversitylibrary.org/bibliography/4360 --- Continued
  as: Memoires Presentés par divers Savants à l'Academie Royale des Sciences
  de l'Institut de France --- Text in French.

\bibitem{DautrayLions1985}
I.~Sneddon, R.~Dautray, M.~Artola, M.~Authier, J.~Lions, P.~Benilan,
  M.~Cessenat, J.~Combes, H.~Lanchon, B.~Mercier, et~al.
\newblock {\em Mathematical Analysis and Numerical Methods for Science and
  Technology: Volume 2 Functional and Variational Methods}.
\newblock Mathematical analysis and numerical methods for science and
  technology. Springer Berlin Heidelberg, 1999.

\bibitem{Strang1972}
G.~Strang.
\newblock Variational crimes in the finite element method.
\newblock In A.~K. Aziz, editor, {\em The Mathematical Foundations of the
  Finite Element Method with Applications to Partial Differential Equations},
  pages 689--710. Academic Press, New York, 1972.
\newblock Proceedings of the Symposium, University of Maryland, Baltimore, Md.,
  1972.

\bibitem{Strang1973}
G.~Strang.
\newblock Piecewise polynomials and the finite element method.
\newblock {\em Bulletin of the American Mathematical Society},
  79(6):1128--1137, 1973.

\bibitem{Young1954}
D.~Young.
\newblock Iterative methods for solving partial difference equations of
  elliptic type.
\newblock {\em Transactions of the American Mathematical Society},
  76(1):92--111, 1954.

\bibitem{Young1972}
D.~M. Young.
\newblock Generalizations of property a and consistent orderings.
\newblock {\em SIAM Journal on Numerical Analysis}, 9(3):454--463, 1972.

\bibitem{ZangWongHussaini1982}
T.~A. Zang, Y.~S. Wong, and M.~Y. Hussaini.
\newblock Spectral multigrid methods for elliptic equations.
\newblock {\em Journal of Computational Physics}, 48:485--501, 1982.

\bibitem{ZangWongHussaini1984}
T.~A. Zang, Y.~S. Wong, and M.~Y. Hussaini.
\newblock Spectral multigrid methods for elliptic equations {II}.
\newblock {\em Journal of Computational Physics}, 54:489--507, 1984.

\end{thebibliography}

\appendix

\section{Relations between the tridiagonal solver and previous work}

\begin{cor}
  Lemma \ref{lem:inversionrichardsonequivalence} is a reformulation of
  \cite[Lemma 1]{GanderKwokZhang2018}
  \label{cor:inversionequivalence}
\end{cor}
\begin{proof}
  We immediately see that $P$ in Lemma
  \ref{lem:inversionrichardsonequivalence} is equal to $P_c$ in
  \cite[Eq. (24)]{GanderKwokZhang2018}, we now look for $R_c$
  \begin{align*}
    R\left(I-M S^{-1}\right) &=
    \left(
    \begin{array}{cc}
      -A^{-1}B & I
    \end{array}
    \right)
    \left(
    \left(
    \begin{array}{cc}
      I & 0 \\
      0 & I
    \end{array}
    \right)
    -
    \left(
    \begin{array}{cc}
      A & B \\
      C & D
    \end{array}
    \right)
    \left(
    \begin{array}{cc}
      A^{-1} & 0 \\
      0 & 0
    \end{array}
    \right)
    \right) \\
    &=
    \left(
    \begin{array}{cc}
      -A^{-1}B & I
    \end{array}
    \right)
    \left(
    \left(
    \begin{array}{cc}
      I & 0 \\
      0 & I
    \end{array}
    \right)
    -
    \left(
    \begin{array}{cc}
      I & 0 \\
      C A^{-1} & 0
    \end{array}
    \right)
    \right)\\
    &=
    \left(
    \begin{array}{cc}
      -A^{-1}B & I
    \end{array}
    \right)
    \left(
    \begin{array}{cc}
      0 & 0 \\
      -C A^{-1} & I
    \end{array}
    \right)\\
    &=
    \left(
    \begin{array}{cc}
      -C A^{-1} & I
    \end{array}
    \right),
  \end{align*}
  and we immediately identify $R_c$.

  It is left to show that $S^{-1} = P_f (R_f M P_f)^{-1} R_f$,
  \begin{align*}
    P_f (R_f M P_f)^{-1} R_f &=
    \left(
    \begin{array}{c}
      I \\
      0
    \end{array}
    \right)
    \left(
    \left(
    \begin{array}{cc}
      I &
      0
    \end{array}
    \right)
    \left(
    \begin{array}{cc}
      A & B \\
      C & D
    \end{array}
    \right)
    \left(
    \begin{array}{c}
      I \\
      0
    \end{array}
    \right)
    \right)^{-1}
    \left(
    \begin{array}{cc}
      I &
      0
    \end{array}
    \right)\\
    &=
    \left(
    \begin{array}{c}
      I \\
      0
    \end{array}
    \right)
    A^{-1}
    \left(
    \begin{array}{cc}
      I &
      0
    \end{array}
    \right) = S^{-1},
  \end{align*}
  and the proof is concluded.
\end{proof}

We show below that cyclic reduction as in \cite{BuzbeeGolubNielson70} is a special case of
the inverse factorization in \cite[Lemma 1]{GanderKwokZhang2018}.
To that end, we keep the notation in \cite{BuzbeeGolubNielson70}.
\begin{lem}[Cyclic reduction as a special case of inverse factorization.]
  Consider the system of equations
  \[
  M x = y
  \]
  where $M$ is an $N \times N$ real symmetric matrix of block tridiagonal form
  \begin{equation*}
    M = 
    \left(\begin{matrix}
        A & T & \dots & 0 & 0 \\
        T & A & \dots & 0 & 0 \\
        \dots & \dots &\dots & \dots & \dots \\
        0 & 0 & \dots & A & T\\ 
        0 & 0 & \dots & T & A\\ 
      \end{matrix}\right).
  \end{equation*}
  The matrices $A$ and $T$ are $p \times p$ symmetric matrices, and we assume
  that
  \[
  AT=TA,
  \]
  the inverse factorization below, as in \cite[Lemma 1]{GanderKwokZhang2018},
  \begin{equation*}
    M^{-1} = P_\text{c}(R_\text{c} M P_\text{c})^{-1} R_\text{c} + P_\text{f}(R_\text{f} M P_\text{f})^{-1} R_\text{f},
  \end{equation*}
  is equivalent cyclic reduction \cite{BuzbeeGolubNielson70},
  \begin{multline}
    \left(\begin{matrix}
        (2T^2-A^2) & T^2 & 0 & \dots & 0 & 0 & 0 \\
        T^2 & (2T^2-A^2) & T^2 & \dots & 0 & 0 & 0\\
        \dots & \dots & \dots & \dots & \dots & \dots & \dots \\
        0 & 0 & 0 & \dots & T^2 & (2T^2-A^2) & T^2 \\
        0 & 0 & 0 & \dots & 0 & T^2 & (2T^2-A^2) 
    \end{matrix}\right)
    \left(\begin{matrix}
        x_2 \\
        x_4 \\
        \dots \\
        x_{N-3} \\
        x_{N-1}
      \end{matrix}\right) \\
    =
    \left(\begin{matrix}
        T y_1 + T y_3 - A y_2 \\
        T y_3 + T y_5 - A y_4 \\
        \dots \\
        T y_{N-6} + T y_{N-4} - A y_{N-3} \\
        T y_{N-2} + T y_{N} - A y_{N-1} \\
      \end{matrix}\right),
    \label{eqn:CyclicReduction1}
  \end{multline}
  \begin{align}
    \left(\begin{matrix}
        A & 0 & 0 & \dots & 0 & 0 & 0 \\
        0 & A & 0 & \dots & 0 & 0 & 0\\
        \dots & \dots & \dots & \dots & \dots & \dots & \dots \\
        0 & 0 & 0 & \dots & 0 & A & 0 \\
        0 & 0 & 0 & \dots & 0 & 0 & A 
    \end{matrix}\right)
    \left(\begin{matrix}
        x_1 \\
        x_3 \\
        \dots \\
        x_{N-2} \\
        x_{N}
      \end{matrix}\right) 
    =&
    \left(\begin{matrix}
        y_1 - T x_2 \\
        y_3 - T y_2 - T y_4 \\
        \dots \\
        y_{N-2} - T y_{N-3} - T y_{N-1} \\
        y_N - T y_{N-1} \\
      \end{matrix}\right).
    \label{eqn:CyclicReduction2}
  \end{align}
  \label{lem:cyclicequivalence}
\end{lem}
\begin{proof}
  Take \cite[Lemma 1]{GanderKwokZhang2018} to solve $M x = y$. Consider
  \begin{align*}
    x =& P_\text{c}(R_\text{c} M P_\text{c})^{-1} R_\text{c} y + P_\text{f}(R_\text{f} M P_\text{f})^{-1} R_\text{f} y,
  \end{align*}
  and verify that $P_\text{f}(R_\text{f} A P_\text{f})^{-1} R_\text{f} y$ does not contribute to even nodes, hence we have
  \begin{align*}
    x_\text{even}=& (R_\text{c} M P_\text{c})^{-1} R_\text{c} y \\
    (R_\text{c} M P_\text{c}) x_\text{even}=& R_\text{c} y,
  \end{align*}
  which is exactly the system in equation \eqref{eqn:CyclicReduction1}.
  Then we have
  \begin{align*}
    x =& P_\text{c}(R_\text{c} M P_\text{c})^{-1} R_\text{c} y + P_\text{f}(R_\text{f} M P_\text{f})^{-1} R_\text{f} y \\
    x =& P_\text{c} x_\text{even} + P_\text{f}(R_\text{f} M P_\text{f})^{-1} R_\text{f} y \\
    x_\text{odd} =& R_\text{f} P_\text{c} x_\text{even} + (R_\text{f} M P_\text{f})^{-1} R_\text{f} y \\   
    (R_\text{f} M P_\text{f}) x_\text{odd} =& R_\text{f} y + (R_\text{f} M P_\text{f}) R_\text{f} P_\text{c} x_\text{even},
  \end{align*}
  which is exactly the system in equation \eqref{eqn:CyclicReduction2}.
\end{proof}

\end{document}